\documentclass[11pt]{amsart}

\usepackage[colorlinks=true, pdfstartview=FitV, linkcolor=blue, 
citecolor=blue]{hyperref}

\usepackage{amssymb,amsmath,amscd}
\usepackage[foot]{amsaddr}
\usepackage{array}
\usepackage{bbm}
\usepackage{xfrac}
\usepackage{graphicx,xcolor}
\usepackage{caption}
\usepackage{subcaption}
\usepackage{a4wide}
\usepackage{lscape}
\usepackage{tikz}
\usepackage{mathtools}
\usepackage{mathrsfs}
\tikzset{every picture/.style={line width=0.75pt}} %set default line width to 0.75pt
\usetikzlibrary{shapes.geometric}
\usetikzlibrary{quotes,angles,positioning}

\newtheorem{theorem}{Theorem}[section]
\newtheorem{prop}[theorem]{Proposition}

\theoremstyle{definition}

\newtheorem{remark}[theorem]{Remark}

\newcommand{\ts}{\hspace{0.5pt}}

\newcommand{\RR}{\mathbb{R}\ts}
\newcommand{\CC}{\mathbb{C}}
\newcommand{\ZZ}{\mathbb{Z}}

\newcommand{\NN}{\mathbb{N}}

\newcommand{\QQ}{\mathbb{Q}}

\newcommand{\XX}{\mathbb{X}}

\newcommand{\ee}{\ts\mathrm{e}}
\newcommand{\ii}{\mathrm{i}\ts}
\newcommand{\bs}{\boldsymbol}

\renewcommand{\phi}{\varphi}
\renewcommand{\theta}{\vartheta}
\renewcommand{\a}{\bs{a}^{{}_{\perp}}}

\newcommand{\exend}{\hfill$\Diamond$}

\newcommand{\myfrac}[2]{\frac{\raisebox{-2pt}{$#1$}}
      {\raisebox{0.5pt}{$#2$}}}

\numberwithin{equation}{section}

\makeatletter
\renewcommand{\@captionfont}{\small}
\makeatother

\DeclareFontFamily{U}{mathx}{\hyphenchar\font45}
\DeclareFontShape{U}{mathx}{m}{n}{ <5> <6> <7> <8> <9> <10>
   <10.95> <12> <14.4> <17.28> <20.74> <24.88> mathx10 }{}
\DeclareSymbolFont{mathx}{U}{mathx}{m}{n}
\DeclareFontSubstitution{U}{mathx}{m}{n}
\DeclareMathAccent{\widecheck}{0}{mathx}{"71}

\newcommand{\defeq}{\mathrel{\mathop:}=}

\begin{document}

\title[Patch frequencies in rhombic Penrose tilings ]
{Patch frequencies in rhombic Penrose tilings }

\author{Jan Maz\'a\v c}
\address{Fakult\"at f\"ur Mathematik, Universit\"at Bielefeld,
  \newline \indent Postfach 100131, 33501 Bielefeld, Germany}
\email{jmazac@math.uni-bielefeld.de}

\begin{abstract} 
This short exposition presents an efficient algorithm for an exact calculation of patch frequencies for the rhombic Penrose tiling. We recall a construction of Penrose tilings via dualisation, and by extending the known method for obtaining vertex configurations, we obtain the desired algorithm. It is then used to determine the frequencies of several particular large patches which appear in the literature. The analogous approach works for a~particular class of tilings and is also explained in detail for the Ammann--Beenker tiling. 
\end{abstract}

\maketitle

\section{Introduction}
The idea of a~non-periodic tiling of a plane with five-fold symmetry goes back to Kepler's famous Figure \texttt{Aa} in \cite{Kepler}. The (rhombic) tiling introduced by  Roger Penrose in \cite{Pen74} is an aperiodic five-fold symmetric tiling of a plane with two prototiles --- a thick and a thin rhombus. There are many ways to generate this tiling. One can define local matching rules, or one can think of it as an inflation tiling and define inflation rules. A more algebraic approach is due to de Bruijn \cite{dBr81}. It relies on the dualisation of a pentagrid, i.e., the union of five rotated lattices. An overview of the methods can be found, for example, in \cite{TAO}. We are interested in another algebraic, yet different, way. It profits from the geometry of the root lattice $A_4$ and the fact that this lattice is a ``minimal" one with five-fold symmetry. Again, this approach uses dualisation; in this scenario, the duality relation between Voronoi and Delone cells (and their complexes).

Recently, the Penrose tiling was considered as an infinite graph and has been studied using tools from graph theory. 
One can consider its graph-theoretic properties like Hamiltonicity, Eulericity, or (perfect) matchings \cite{FSP20,LBSPF22}, but one can also assign an operator acting on this graph and study its spectral properties. In \cite{DEFM22}
Damanik et al. study the properties of a Laplacian on various tilings, among them the rhombic Penrose one. 
They studied a~tile model for the Laplacian, and they were able to show some examples of locally-supported eigenfunctions, which are also known from other papers \cite{FATK88}. Recently, Oktel published several papers dealing with a similar problem for the vertex model for different tilings \cite{O21,O22a,O22b}. 
For all these models, one can further study the \textit{integrated density of states} (IDS), which is a function that counts the number of states (different eigenfunctions) up to a given energy. 
It was shown that this function is discontinuous. More precisely, if one can find a locally supported eigenfunction with energy $E$ of the Laplacian, the IDS has  a~discontinuity jump at the point $E$. 
The size of this gap is at least as big as the frequency of the eigenfunction's support, i.e., the frequency of the corresponding patch. \cite{DEFM22}
Thus, knowing the frequency, one gets a lower bound on the size of the gap. Damanik et al. used a direct approach to calculate the frequencies. 
Namely, they count the number of occurrences of the support of a given eigenfunction in growing approximants of the entire tiling. 
The same method was employed earlier by Fujiwara et al. \cite{FATK88}. 
There is an obvious disadvantage of this method. 
Indeed, one has to deal with the boundary of the approximants, which may include parts of the studied patch. 
Another problem constitutes the way of choosing the approximants.
Lastly, the resulting frequency is always given as a numerical approximation.
Therefore, we aim to fill this gap by showing an algebraic way to obtain the frequencies of arbitrary finite patches in (not only) Penrose rhombic tiling exactly, without any need for the inflation method. For Penrose rhombic tilings, there already exists a method by Zobetz and Preislinger \cite{ZP90} using de Bruijn's approach, which enables a calculation of frequencies of vertex configurations in generalized Penrose tilings. Still, our approach provides a more general framework as it allows us to effectively calculate an exact frequency of arbitrary large patches for a wider class of tilings. As far as we are aware, there does not exist any algorithm that would actually enable the calculation of exact frequencies for arbitrary finite patches. 

This paper is structured as follows. In Section \ref{sec:A4lattice}, we recall the geometry of the $A_4$ lattice and its Voronoi complex and of their dual objects. Further, in Section \ref{sec:fivefold}, we recall a representation of a cyclic group of order 5 (which acts naturally on the lattice $A_4$), which exhibits five-fold symmetry in a plane. Section \ref{sec:dual} evokes the dualisation method and its benefits. These sections are almost fully based on \cite{BKSZ90}. We recall them as they are necessary for the algorithm. The crucial point is that it describes tilings rather than point sets by a variant of the projection method known as dualisation. In particular, the standard model set approach via the intersection of translated windows \cite[Cor. 7.3]{TAO} is practically unable to give the frequencies of large patches. The algorithm for determining the frequencies is presented in Section \ref{sec:freq}. In Appendix 1, we apply it to several patches coming from \cite{DEFM22}. The second appendix is devoted to a brief summary of the patch frequencies in Ammann--Beenker tilings. 

\section{The root lattice $A_4$, its dual, and their properties}
\label{sec:A4lattice}

\begin{figure}
    \centering
 
\tikzset{every picture/.style={line width=0.75pt}} %set default line width to 
\begin{tikzpicture}[x=0.75pt,y=0.75pt,yscale=-1,xscale=1]
\filldraw [black] (100,160) circle (4pt);
\filldraw [black] (160,160) circle (4pt);
\filldraw [black] (220,160) circle (4pt);
\filldraw [black] (280,160) circle (4pt);

\draw    (100,160) -- (280,160) ;

\draw (100,180) node [anchor=north][inner sep=0.75pt]    {\small $\bs{e}_1-\bs{e}_2$};
\draw (160,180) node [anchor=north][inner sep=0.75pt]    {\small$\bs{e}_2-\bs{e}_3$};
\draw (220,180) node [anchor=north][inner sep=0.75pt]    {\small$\bs{e}_3-\bs{e}_4$};
\draw (280,180) node [anchor=north][inner sep=0.75pt]    {\small$\bs{e}_4-\bs{e}_5$};
\end{tikzpicture}

\caption{The Dynkin diagram $A_4$. Every node represents a basis vector, and their geometry is encoded via the lines. If two vertices are connected, their scalar product is -1. Otherwise, they are orthogonal.}
\label{fig:dynkin}
\end{figure}
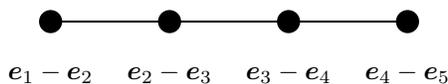

The lattice $A_4$ can be understood in different ways. Perhaps the most natural one (explaining its name) is that $A_4$ is the root lattice of the semisimple Lie algebra $\mathfrak{sl}_5(\CC)$. 
On the other hand, its explicit description as an intersection of the primitive 5-dimensional cubic lattice with a~4-dimensional hyperplane allows us to simplify some calculations. 
Thus, let $\bs{e}_1,\dots,\bs{e}_5$ be the standard basis vectors of $\RR^5$ and set $\bs{s} = \bs{e}_1+\cdots + \bs{e}_5$. 
Let further $\mathcal{S}=\{\bs{x}\in\RR^5 \ : \ \bs{s}\cdot \bs{x} = 0  \} \simeq \RR^4$, be a~4-dimensional hyperplane in $\RR^5$. Then, one has
\[ A_4\, =\, \mathcal{S}\cap \ZZ^5. \]
The resulting lattice is generated by four vectors, namely
\[ A_4 = \langle \bs{e}_1-\bs{e}_2,\, \bs{e}_2-\bs{e}_3,\, \bs{e}_3-\bs{e}_4,\, \bs{e}_4-\bs{e}_5 \rangle_{\ZZ}.  \]
Alternatively, we can depict the root lattice $A_4$ as a Dynkin diagram; see Figure \ref{fig:dynkin}. 

Note that the generating vectors $\bs{e}_i-\bs{e}_{i+1}$ are fundamental (or simple) roots of the root system of $\mathfrak{sl}_5(\CC)$. This system consists of 20 root vectors, namely $\bs{e}_i - \bs{e}_j$ with $1\leqslant i,j \leqslant 5 $ and $i\neq j$. 
For our further analysis, we need to describe the maximal point symmetry group $H_{A_4}$ at the origin of the lattice $A_4$. It is isomorphic with the automorphism group of the generating root system. 
The root system is, by definition, invariant under the action of the Weyl group $W(A_4)$, which is the permutation group $S_5$ in this case. 
Moreover, central inversion is an additional symmetry generating the group $Z_2$. 
Thus the group $H_{A_4}$ is isomorphic to 
\[ H_{A_4} \, \simeq \, W(A_4) \times Z_2 \, \simeq \, S_5 \times Z_2. \]

The 20 root vectors also determine the \emph{Voronoi cell} $V_{A_4}(\bs{0})$ around the origin, i.e., all vectors in the underlying hyperplane $\mathcal{S}$ which are not further apart (with respect to the Euclidean distance) from the origin than to any other lattice point, so
\[V_{A_4}(\bs{0}) = \{ \bs{x} \in \mathcal{S} \ : \ \forall \bs{v} \in A_4  \ : \ ||\bs{x}-\bs{v}|| \geqslant ||\bs{v}||  \}.  \]
The Voronoi cell can also be understood as an intersection of closed half-spaces $H^+_{\bs{v}}$ corresponding to $\bs{v} \in A_4$ defined as   $H^+_{\bs{v}} \defeq  \{ \bs{x} \in \mathcal{S} \ : \ ||\bs{x}-\bs{v}|| \geqslant ||\bs{v}||  \} $. Here, the  Voronoi cell $V_{A_4}(\bs{0})$ is fully determined by the 20 root vectors, i.e., one has 
\[ V_{A_4}(\bs{0}) \, = \, \bigcap_{i\neq j} H^+_{\bs{e}_i-\bs{e}_j }. \]

To obtain a more explicit description of the Voronoi cell $V_{A_4}(\bs{0})$, we have to employ the dual lattice $A^*_4$ and its fundamental domain. The dual lattice can be obtained in many ways. Following Conway's approach via \emph{glue vectors} \cite{Conway_Sloane},
one has
\[A^*_4 \, = \, \bigcup_{i=0}^4 \ [i]+ A_4\, , \]
with the glue vectors 
\[[0]= \begin{pmatrix}0\\0\\0\\0\\0  \end{pmatrix}, \ 
[1] = \myfrac{1}{5}\begin{pmatrix}1\\1\\1\\1\\-4  \end{pmatrix}, \ 
[2] = \myfrac{1}{5}\begin{pmatrix}2\\2\\2\\-3\\-3  \end{pmatrix}, \
[3] = \myfrac{1}{5}\begin{pmatrix}3\\3\\-2\\-2\\-2  \end{pmatrix}, \
[4] = \myfrac{1}{5}\begin{pmatrix}4\\-1\\-1\\-1\\-1  \end{pmatrix}. \
\]

This description allows one to immediately recognise $A_4$ as a proper sublattice in its dual lattice $A^*_4$. 
Moreover, the quotient group $\sfrac{A^*_4}{A_4}\simeq C_5$ is of order 5, and the representatives can be chosen as the glue vectors. 
On the other hand, for upcoming calculations, it is convenient to write down the generators of the lattice. 
Here, $A^*_4$ is spanned by the vectors 
\[ \bs{a}_i \, = \, \bs{e}_i - \myfrac{1}{5}\bs{s}\]
with $1\leqslant i \leqslant 5$ and $\bs{s}$ as above. Note that the generating vectors are not linearly independent since $\sum_{i=1}^5\bs{a}_i = \bs{0}$. Finally, one can use them to describe the Voronoi cell 
\[V_{A_4}(\bs{0}) = \bigg\{\bs{x} \in \mathcal{S} \ : \ \bs{x} = \myfrac{1}{2} \sum_{i=1}^5 \lambda_i \bs{a}_i\, , \ |\lambda_i| \leqslant 1 \bigg\}.  \]
This object is a regular 4-dimensional convex polytope, sometimes considered as a dual polytope to the \emph{runcinated 5-cell}. It has the full symmetry $W(A_4) \times Z_2$. The polytope possesses 30 vertices, 70 bounding edges, 60 bounding polygons (i.e., polytopes of dimension 2) and 20 bounding polytopes of dimension 3. 
Henceforth, we refer to them as \emph{$k$-boundaries}, with $0\leqslant k \leqslant 3$. 
Baake et al. \cite{BKSZ90} provide a~careful analysis of all $k$-boundaries and their explicit description together with one of their corresponding duals in the sense of \cite{KS89}.  
Important to us here are the 2-boundaries, the vertices, and the corresponding dual objects as follows.

The 2-boundary polygons are given by 
\[ P(++-\bigcirc \bigcirc) = \left\{ \myfrac{1}{2} \bigl( \bs{a}_1+\bs{a}_2 - \bs{a}_3 + \lambda_4\bs{a}_4 + \lambda_5\bs{a}_5 \bigr) : \ |\lambda_i| \leqslant 1    \right\} \]
together with all polygons obtained via vertex permutations and sign flips. There is an explicit action of the group $H_{A_4}$ on the set of 2-boundaries.  This action can be  encoded on the level of the signature $(++-\bigcirc \bigcirc)$ as well. In particular, a permutation just permutes the indices, and a sign flip affects the signs and $\bigcirc$ remains unchanged. 
From the geometric point of view, $P(++-\bigcirc \bigcirc)$ is a rhombus; therefore, it will play a crucial role in constructing the Penrose rhombus tiling. 
The 2-boundary dual to $P(++-\bigcirc \bigcirc)$ is the triangle $P^{*}(++-\bigcirc \bigcirc)$ defined as
\[ P^{*}(++-\bigcirc \bigcirc) = \left\{ \mu_1(\bs{a}_1-\bs{a}_3) +  \mu_2(\bs{a}_2 -\bs{a}_3) : \ \mu_i\geqslant 0, \ \mu_1+\mu_2 \leqslant 1    \right\}.  \]

The correspondence between $P$ and $P^{*}$ is one-to-one, and the boundaries intersect with their duals at precisely one point. 

The 30 vertex points of the Voronoi cell $V_{A_4}(\bs{0})$ are exactly those points of $\mathcal{S}$ with the largest distance to the lattice $A_4$. 
In terms of the theory of root lattices, they are called \emph{holes} \cite{Conway_Sloane}. 
Points with the maximal possible distance to $A_4$ are called \emph{deep holes}, and the remaining ones are \emph{shallow holes}. In our case, the vertex 
\[P(++++-) = \myfrac{1}{2}\bigl(\bs{a}_1 +\bs{a}_2+\bs{a}_3+\bs{a}_4-\bs{a}_5 \bigr) = \bs{a}_1 +\bs{a}_2+\bs{a}_3+\bs{a}_4 \]
and all its images under $W(A_4) \times Z_2$ are  the shallow holes, whereas the 20 points of type 
\[P(+++--) =  \myfrac{1}{2}\bigl(\bs{a}_1 +\bs{a}_2+\bs{a}_3-\bs{a}_4-\bs{a}_5 \bigr) = \bs{a}_1 +\bs{a}_2+\bs{a}_3  \]
are the deep holes.

The dual objects to deep and shallow holes are four-dimensional cells. Namely, one gets a~4-dimensional~simplex
\begin{equation}
\label{eq:4D_dual_1}
P^{*}(++++-) = \left\{ \sum_{i=1}^4\mu_i(\bs{a}_i -\bs{a}_5) : \mu_i \geqslant 0,  \  \sum_{i=1}^4\mu_i \leqslant 1 \right\}, 
\end{equation}
and a 4-dimensional Archimedian polytope 

\begin{equation}
\begin{split}
\label{eq:4D_dual_2}
    P^{*}(+++--) = \left\{ \sum_{i=1}^3\mu_i(\bs{a}_i -\bs{a}_4) \right. &\left. + \sum_{i=1}^3\mu_{i+3}(\bs{a}_i -\bs{a}_5) \ : \ \mu_i \geqslant 0,  \  \sum_{i=1}^3 \mu_i \leqslant 1, \right. \\  &\left. \sum_{i=4}^6\mu_i \leqslant 1, \ \mu_i+\mu_{i+3}\leqslant 1 \mbox{ for all } 1\leqslant i \leqslant 3 \right\},
\end{split}
\end{equation}
and all their images under the symmetry operations $H_{A_4}. $

Since $A_4$ is a lattice, one has the same vertex configuration around any of its points up to translation. Thus, the Voronoi cell $V_{A_4}(\bs{v})$ around $\bs{v}$ is a translate $V_{A_4}(\bs{0}) + \bs{v}$. Further, one can collect all $k$-boundaries and think of them in terms of complexes. In particular, one can define the \emph{Voronoi complex} 
\[\mathscr{V}\defeq \left\{ P \subset \mathcal{S} \ : \ P \mbox{ is a boundary of some }V_{A_4}(\bs{v}) \ \mbox{with} \ \bs{v} \in A_4 \right\} \]
and for $0\leqslant k \leqslant 4$ its \emph{$k$-skeleton} 
\[\mathscr{V}^{(k)}\defeq \left\{ P \in \mathscr{V} \ : \ P \mbox{ is a $k$-boundary}\right\}. \]
The properties of the duality leads to the \emph{dual Voronoi complex} and its \emph{dual $k$-skeleton } as
\begin{align*}
    \mathscr{V}^{*} & \defeq \left\{P^{*} \ : \ P\in\mathscr{V} \right\},\\
    \mathscr{V}^{*(k)} & \defeq \left\{ P^{*} \in \mathscr{V}^{*} : P \mbox{ is a $(4{-}k)$-boundary} \right\}. 
\end{align*}

Taking any vertex $\bs{v}^{*}$ of the Voronoi cell $V_{A_4}(\bs{v})$ for some $\bs{v} \in A_4$, i.e., $\bs{v}^{*} \in \mathscr{V}^{(0)}$, the associated dual object, which is a full 4D polytope, will be denoted by $V^{*}(\bs{v}^{*})$ as it plays a~similar role as the Voronoi cell. 

As mentioned above, different points appear within the point sets studied. We have to deal with points of the lattice $A_4$ and with the vertices of its Voronoi cells. The latter split into two categories, deep and shallow holes. In order to distinguish them, one can introduce a modulo function $r$ defined for any  point $\bs{v}^{*} = \sum_{i=1}^4 n_i \bs{a}_i \in A_4^{*}$ as 
\[ r(\bs{v}^{*})\defeq \left( \sum_{i=1}^4 n_i \right) \mod 5.  \]
It is clear that $r: A^{*}_4 \longrightarrow \sfrac{\ZZ}{5\ZZ}$. Since the generating vectors $\bs{e}_i - \bs{e}_{i+1}$ of the lattice $A_4$ fulfil 
\[ \bs{e}_i - \bs{e}_{i+1} = \bs{a}_i - \bs{a}_{i+1}, \]
one has immediately 
\[ A_4 = \ker (r). \]
Further, one obtains the characterisation of shallow and deep holes in terms of $r(\bs{v}^{*})$. In particular, 
\begin{align*}
    \bs{v}^{*} \mbox{is a shallow hole } & \Leftrightarrow \  r(\bs{v}^{*}) =\pm 1 \mod 5, \\
    \bs{v}^{*} \mbox{is a deep hole } & \Leftrightarrow \  r(\bs{v}^{*}) =\pm 2 \mod 5.
\end{align*} 
\begin{remark}
The function $r$ corresponds to the \textit{index} function in de Bruijn's construction \cite{dBr81}. This is not surprising because de Bruijn's construction implicitly uses root lattice as a~Minkowski embedding of fifth roots of unity as explained in \cite[Sec. 7.5.2]{TAO}. \exend
\end{remark}

\section{Representation with five-fold symmetry }
\label{sec:fivefold}
We have already mentioned that $W(A_4)$ acts on the generators of $A_4$ via  permutations of the basis vectors $\bs{e}_i$. This action has two invariant subspaces, namely $\CC\bs{s}$ and $\mathcal{S}$. 
The linear representation of $S_5 \simeq W(A_4)$ is irreducible on $\mathcal{S}$, and to find a real irreducible representation capturing the fivefold symmetry in plane, one has to restrict oneself to a suitable subgroup. 
Therefore, consider the cyclic group $C_5$, a subgroup of $W(A_4)$. Its generating element $g=(12345)$ acts on the basis $(\bs{e}_1,\dots,\bs{e}_5)$ via  the matrix
\[D(g) = \begin{pmatrix} 
0& 0 & 0& 0& 1 \\
1& 0 & 0& 0& 0 \\
0& 1 & 0& 0& 0 \\
0& 0 & 1& 0& 0 \\
0& 0 & 0& 1& 0 \\
\end{pmatrix}.\]
To find the possible representations means to find the real Jordan form of $D(g)$ via an orthogonal matrix J. The real Jordan form reads
\[ JD(g)J^{-1} = \begin{pmatrix} 
\cos\frac{2\pi}{5} & -\sin\frac{2\pi}{5} & 0 & 0 & 0\\[5pt]
\sin\frac{2\pi}{5} & \cos\frac{2\pi}{5} &0 & 0 & 0 \\[5pt]
0 & 0 & \cos\frac{4\pi}{5} & -\sin\frac{4\pi}{5} & 0 \\[5pt]
0 & 0 & \sin\frac{2\pi}{5} & \cos\frac{2\pi}{5} & 0 \\[5pt]
0&0&0&0&1
\end{pmatrix} = D^{\parallel}(g) \oplus D^{\perp}(g) \oplus D^0(g)\]
and provides three irreducible real representations $D^{\parallel}(g),\ D^{\perp}(g), \ D^0(g)$. 
The matrix $J$ read columnwise provides a new basis as one can directly read from 
\[ JD(g) = \bigl(D^{\parallel}(g) \oplus D^{\perp}(g) \oplus D^0(g)\bigr)J. \]
In particular, one has
\[J = \sqrt{\frac{2}{5}} \begin{pmatrix}
1 & \cos\frac{2\pi}{5} & \cos\frac{4\pi}{5} & \cos\frac{4\pi}{5} & \cos\frac{2\pi}{5} \\[5pt]
0 & \sin\frac{2\pi}{5} &\sin\frac{4\pi}{5} & -\sin\frac{4\pi}{5} & -\sin\frac{2\pi}{5}\\[5pt]
1 &\cos\frac{4\pi}{5} & \cos\frac{2\pi}{5} & \cos\frac{2\pi}{5} &\cos\frac{4\pi}{5} \\[5pt]
0 &\sin\frac{4\pi}{5} & -\sin\frac{2\pi}{5}& \sin\frac{2\pi}{5}& -\sin\frac{4\pi}{5} \\[5pt]
\sqrt{\frac{1}{2}} & \sqrt{\frac{1}{2}} &\sqrt{\frac{1}{2}} & \sqrt{\frac{1}{2}}& \sqrt{\frac{1}{2}}
\end{pmatrix}.\] 

Since the trivial representation $D^{(0)}(g)$ is carried by the subspace $\CC\bs{s}$, it follows that $D^{\parallel}(g)$ and $ D^{\perp}(g)$ are contained in $\mathcal{S}$. Thus, one has to decompose $\mathcal{S}$ as a direct sum of two subspaces, $\mathcal{S}^{\parallel}$ and $\mathcal{S}^{\perp}$. The representation of $g$ in $\mathcal{S}^{\parallel}$ and $\mathcal{S}^{\perp}$ is a rotation about $\tfrac{2\pi}{5}$, and $\tfrac{4\pi}{5}$, respectively. 

Denote by $\pi_{{}_\parallel}$ and $\pi_{{}_\perp}$ the projections from $\mathcal{S}$ onto $\mathcal{S}^{\parallel}$ and $\mathcal{S}^{\perp}$, respectively. The projection of basis vectors $\pi_{{}_\parallel}(\bs{e}_i)$ are given by the first and second row of $i$-th column of $J$, and $\pi_{{}_\perp}(\bs{e}_i)$ are given by the third and fourth row of the same column. Figure \ref{fig:basis_projections} depicts the projections of the basis vectors $\bs{e}_i$, which exhibit the desired five-fold symmetry. Since $\sum_{i=1}^5 \pi_{{}_\parallel}(\bs{e}_i) = \sum_{i=1}^5 \pi_{{}_\perp}(\bs{e}_i) = \bs{0}$, one immediately gets $\pi_{{}_\parallel}(\bs{a}_i) = \pi_{{}_\parallel}(\bs{e}_i)$ and $\pi_{{}_\perp}(\bs{a}_i) = \pi_{{}_\perp}(\bs{e}_i)$ for all $1\leqslant i \leqslant 5$. 

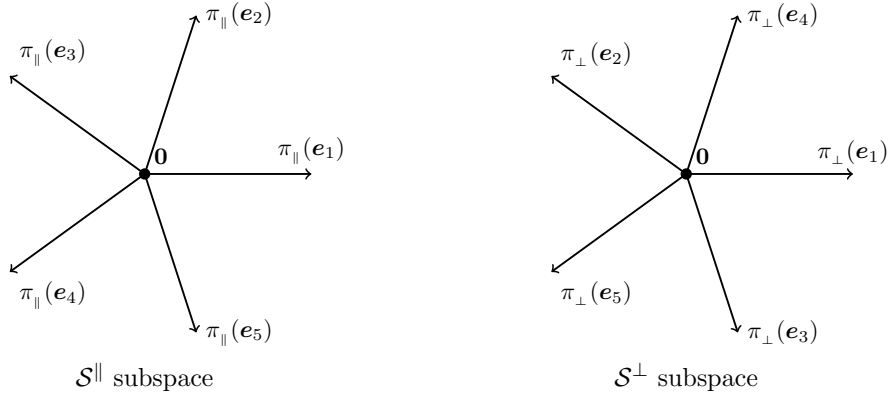
\begin{figure}
    \centering
    \scalebox{0.9}{
\begin{tikzpicture}
			\filldraw[black] (0,0) circle (2pt) node[above right]{$\small\bs{0}$}; 
			\node[regular polygon, regular polygon sides=5, shape border rotate = -18,
			inner sep=1.4cm] (s) at (0,0) {} ;
			\draw[thick,  ->] (s.center) -- (s.corner 5)
			node [above]{$\small \pi_{{}_\parallel}(\bs{e}_1)$};
			\draw[thick, ->] (s.center) -- (s.corner 1)
			node [right]{ $\small \pi_{{}_\parallel}(\bs{e}_2)$};
			\draw[thick,  ->] (s.center) -- (s.corner 2)
			node [above right]{$\small \pi_{{}_\parallel}(\bs{e}_3)$};
			\draw[thick, ->] (s.center) -- (s.corner 3)
			node [below right]{$\small \pi_{{}_\parallel}(\bs{e}_4)$};
			\draw[thick, ->] (s.center) -- (s.corner 4)
			node [right]{$\small \pi_{{}_\parallel}(\bs{e}_5)$};
			\node (A) at (0,-3) {$\mathcal{S}^{\parallel}$ subspace};
			\filldraw[black] (8,0) circle (2pt) node[above right]{$\small\bs{0}$}; 
			\node[regular polygon, regular polygon sides=5, shape border rotate = -18,
			inner sep=1.4cm] (r) at (8,0) {} ;
			\draw[thick,  ->] (r.center) -- (r.corner 5)
			node [above]{$\small \pi_{{}_\perp}(\bs{e}_1)$};
			\draw[thick, ->] (r.center) -- (r.corner 1)
			node [right]{ $\small \pi_{{}_\perp}(\bs{e}_4)$};
			\draw[thick,  ->] (r.center) -- (r.corner 2)
			node [above right]{$\small \pi_{{}_\perp}(\bs{e}_2)$};
			\draw[thick, ->] (r.center) -- (r.corner 3)
			node [below right]{$\small \pi_{{}_\perp}(\bs{e}_5)$};
			\draw[thick, ->] (r.center) -- (r.corner 4)
			node [right]{$\small \pi_{{}_\perp}(\bs{e}_3)$};
			\node (B) at (8,-3) {$\mathcal{S}^{\perp}$ subspace};
\end{tikzpicture}
}
\caption{Projections of the standard basis $\bs{e}_1, \dots, \bs{e}_5$ into the two subspaces $\mathcal{S}^{\parallel}$ and $\mathcal{S}^{\perp}$, respectively.  }
    \label{fig:basis_projections}
\end{figure}
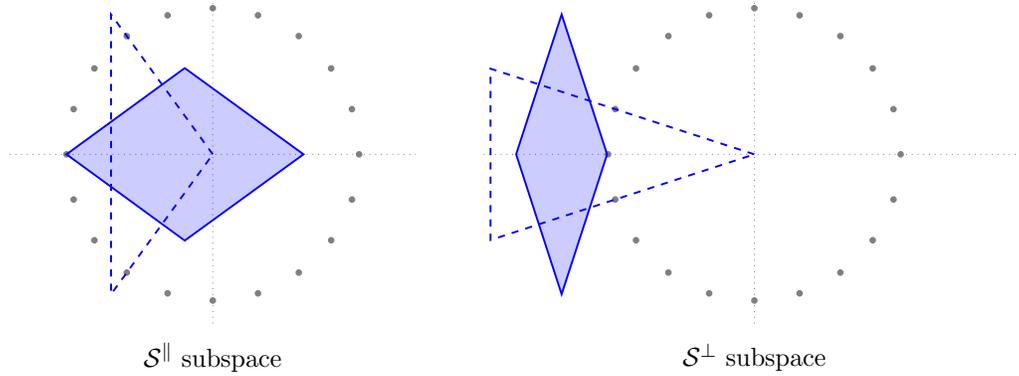
\begin{figure}[h]
    \centering
    \scalebox{0.9}{
\begin{tikzpicture}
    \draw[thin, gray,dotted] (-3,0)--(3,0);
			\draw[thin, gray,dotted] (0,-2.5)--(0,2.5);
			\node[regular polygon, regular polygon sides=20, shape border rotate = -81, inner sep=1.5cm] (s) at (0,0) {} ;
			\node[regular polygon, regular polygon sides=20, shape border rotate = -81, inner sep=0.927cm] (r) at (0,0) {} ;  
			\foreach \x in {1,...,20}
			\node[circle,draw=gray, fill=gray, inner sep=0pt,minimum size=2pt] (s.\x) at (s.corner \x) {};
			\draw[thick,blue, fill=blue, fill opacity = 0.2] (r.corner 1) -- (r.corner 7)-- (s.corner 11)-- (r.corner 15) -- cycle;
			\draw[thick, blue, dashed] (s.center) -- (-1.5029,2.0686)-- (-1.5029,-2.0686) -- cycle;
			\draw[thin, gray,dotted] (4,0)--(12,0);
			\draw[thin, gray,dotted] (8,-2.5)--(8,2.5);
			\node[regular polygon, regular polygon sides=20, shape border rotate = -81, inner sep=1.5cm] (S) at (8,0) {} ;
			\node[regular polygon, regular polygon sides=20, shape border rotate = -81, inner sep=2.853cm] (b) at (8,0) {} ;
			\foreach \x in {1,...,20}
			\node[circle,draw=gray, fill=gray, inner sep=0pt,minimum size=2pt] (S.\x) at (S.corner \x) {};
			\draw[thick,blue, fill=blue, fill opacity = 0.2] (5.825,0) -- (5.1529,2.0684) -- (4.4808,0) --(5.1529,-2.0684) -- cycle;
			\draw[thick, blue, dashed] (b.center) -- (b.corner 10) -- (b.corner 12) -- cycle;
			\node (A) at (0,-3) {$\mathcal{S}^{\parallel}$ subspace};
			\node (B) at (8,-3) {$\mathcal{S}^{\perp}$ subspace};
\end{tikzpicture}
}
    \caption{Images of the 2-boundary $P(-+\bigcirc \bigcirc+)$ and its dual $P^{*}(-+\bigcirc \bigcirc+)$ under the projections $\pi_{{}_\parallel}$ and $\pi_{{}_\perp}$, respectively. The solid blue rhombuses correspond to projections of the 2-boundary, whereas the dashed line indicates the projection of its dual. The gray points are the $20^{\mathrm{th}}$ roots of unity scaled by  $\sqrt{\tfrac{2}{5}}$.   }
    \label{fig:projections_boundaries}
\end{figure}
Projecting the 2-boundary $P$ and its dual 2-boundary $P^{*}$ in both spaces results in a set of triangles and rhombuses, which we use later for the construction of the Penrose tiling. Figure \ref{fig:projections_boundaries} shows the projections of $P(-+\bigcirc \bigcirc+)$ and $P^{*}(-+\bigcirc \bigcirc+)$.
Note that the rhombus vertices always consist of one projection of a shallow hole and three projections of deep holes. 
The position of the shallow hole will later be needed to distinguish different patterns.

\section{Dualisation method}
\label{sec:dual}
One can obtain a space tiling via the so-called \emph{dualisation method}. This method was described in detail in \cite{KS89}, and \cite{TAO} provides an illustrative overview. 
To employ this method, one needs a~Voronoi complex $\mathscr{V}$, its dual (Delone) complex $\mathscr{V}^{*}$, and a suitable \emph{cutting plane}, which carries the desired tiling. 
To get a non-periodic tiling, one has to choose the cutting plane so that it contains at most one lattice point. 

\begin{figure}[h]
    \centering
    \scalebox{0.7}{
\begin{tikzpicture}
    \draw[thin, gray,dotted] (-3,0)--(3,0);
    \draw[thin, gray,dotted] (0,-3)--(0,3);
	\node[regular polygon, regular polygon sides=20, shape border rotate = -81, inner sep=1.5cm] (s) at (0,0) {} ;
	\node[regular polygon, regular polygon sides=20, shape border rotate = -81, inner sep=0.927cm] (r) at (0,0) {} ;
		\foreach \x in {1,...,20}
			    \node[circle,draw=gray, fill=gray, inner sep=0pt,minimum size=2pt] (s.\x) at (s.corner \x) {};
		\draw[thick,blue, fill=blue, fill opacity = 0.2] (r.corner 1) -- (r.corner 7)-- (s.corner 11)-- (r.corner 15) -- cycle;
		\draw[thick, blue, dashed] (s.corner 1) -- (r.corner 5)-- (r.corner 11)-- (r.corner 17) -- cycle;
			\node[circle,draw=red, fill=red, inner sep=1pt, minimum size= 4pt] (Sh) at (s.corner 11) {};
			\node[circle,draw=red, fill=red, inner sep=1pt, minimum size= 4pt] (Sh) at (s.corner 1) {};   
    \draw[thin, gray,dotted] (4.5,0)--(10.5,0);
    \draw[thin, gray,dotted] (7.5,-3)--(7.5,3);
	\node[regular polygon, regular polygon sides=20, shape border rotate = -81, inner sep=1.5cm] (s) at (7.5,0) {} ;
	\node[regular polygon, regular polygon sides=20, shape border rotate = -81, inner sep=0.927cm] (r) at (7.5,0) {} ;
		\foreach \x in {1,...,20}
			    \node[circle,draw=gray, fill=gray, inner sep=0pt,minimum size=2pt] (s.\x) at (s.corner \x) {};
		\draw[thick,blue, fill=blue, fill opacity = 0.2] (r.corner 5) -- (r.corner 11)-- (s.corner 15)-- (r.corner 19) -- cycle;
		\draw[thick, blue, dashed] (s.corner 5) -- (r.corner 9)-- (r.corner 15)-- (r.corner 1) -- cycle;
			\node[circle,draw=red, fill=red, inner sep=1pt, minimum size= 4pt] (Sh15) at (s.corner 15) {};
			\node[circle,draw=red, fill=red, inner sep=1pt, minimum size= 4pt] (Sh5) at (s.corner 5) {};
	\draw[thin, gray,dotted] (12,0)--(18,0);
    \draw[thin, gray,dotted] (15,-3)--(15,3);
	\node[regular polygon, regular polygon sides=20, shape border rotate = -81, inner sep=1.5cm] (s) at (15,0) {} ;
	\node[regular polygon, regular polygon sides=20, shape border rotate = -81, inner sep=0.927cm] (r) at (15,0) {} ;
		\foreach \x in {1,...,20}
			    \node[circle,draw=gray, fill=gray, inner sep=0pt,minimum size=2pt] (s.\x) at (s.corner \x) {};
		\draw[thick,blue, fill=blue, fill opacity = 0.2] (r.corner 9) -- (r.corner 15)-- (s.corner 19)-- (r.corner 3) -- cycle;
		\draw[thick, blue, dashed] (s.corner 9) -- (r.corner 13)-- (r.corner 19)-- (r.corner 5) -- cycle;
			\node[circle,draw=red, fill=red, inner sep=1pt, minimum size= 4pt] (Sh19) at (s.corner 19) {};
			\node[circle,draw=red, fill=red, inner sep=1pt, minimum size= 4pt] (Sh9) at (s.corner 9) {};  
	\node[align=center] (A) at (0,-3.5) {$P(-+\bigcirc \bigcirc+)$ };	
        \node[align=center] (B) at (7.5,-3.5) {$P(+-+\bigcirc \bigcirc)$ };
		\node[align=center] (C) at (15,-3.5) {$P(\bigcirc +-+ \bigcirc)$ };
		
	\draw[thin, gray,dotted] (0.75,-7)--(6.75,-7);
    \draw[thin, gray,dotted] (3.75,-10)--(3.75,-4);
	\node[regular polygon, regular polygon sides=20, shape border rotate = -81, inner sep=1.5cm] (s) at (3.75,-7) {} ;
	\node[regular polygon, regular polygon sides=20, shape border rotate = -81, inner sep=0.927cm] (r) at (3.75,-7) {} ;
		\foreach \x in {1,...,20}
			    \node[circle,draw=gray, fill=gray, inner sep=0pt,minimum size=2pt] (s.\x) at (s.corner \x) {};
		\draw[thick,blue, fill=blue, fill opacity = 0.2] (r.corner 7) -- (r.corner 13)-- (r.corner 19)-- (s.corner 3) -- cycle;
		\draw[thick, blue, dashed] (s.corner 13) -- (r.corner 17)-- (r.corner 3)-- (r.corner 9) -- cycle;
			\node[circle,draw=red, fill=red, inner sep=1pt, minimum size= 4pt] (Sh19) at (s.corner 3) {};
			\node[circle,draw=red, fill=red, inner sep=1pt, minimum size= 4pt] (Sh9) at (s.corner 13) {};  
	\draw[thin, gray,dotted] (8.25,-7)--(14.25,-7);
    \draw[thin, gray,dotted] (11.25,-10)--(11.25,-4);
	\node[regular polygon, regular polygon sides=20, shape border rotate = -81, inner sep=1.5cm] (s) at (11.25,-7) {} ;
	\node[regular polygon, regular polygon sides=20, shape border rotate = -81, inner sep=0.927cm] (r) at (11.25,-7) {} ;
		\foreach \x in {1,...,20}
			    \node[circle,draw=gray, fill=gray, inner sep=0pt,minimum size=2pt] (s.\x) at (s.corner \x) {};
		\draw[thick,blue, fill=blue, fill opacity = 0.2] (r.corner 11) -- (r.corner 17)-- (r.corner 3)-- (s.corner 7) -- cycle;
		\draw[thick, blue, dashed] (s.corner 17) -- (r.corner 1)-- (r.corner 7)-- (r.corner 13) -- cycle;
			\node[circle,draw=red, fill=red, inner sep=1pt, minimum size= 4pt] (Sh19) at (s.corner 7) {};
			\node[circle,draw=red, fill=red, inner sep=1pt, minimum size= 4pt] (Sh9) at (s.corner 17) {};  
	\node[align=center] (C) at (3.75,-10.5) {$P(\bigcirc \bigcirc+-+)$ };	
    \node[align=center] (D) at (11.25,-10.5) {$P(+\bigcirc \bigcirc+-)$ };

\end{tikzpicture}
}  
    \caption{Projections of the different (modulo translation) 2-boundaries $P$ in the $\mathcal{S}^{\parallel}$ which result in a thick rhombus. The solid rhombi correspond to the label, whereas the dashed rhombi are their space inversion. The red point attached to a given rhombus indicates the shallow hole. The gray points are the $20^{\mathrm{th}}$ roots of unity scaled by the factor $\sqrt{\tfrac{2}{5}}$.  }
    \label{fig:projections_rhombi}
\end{figure}
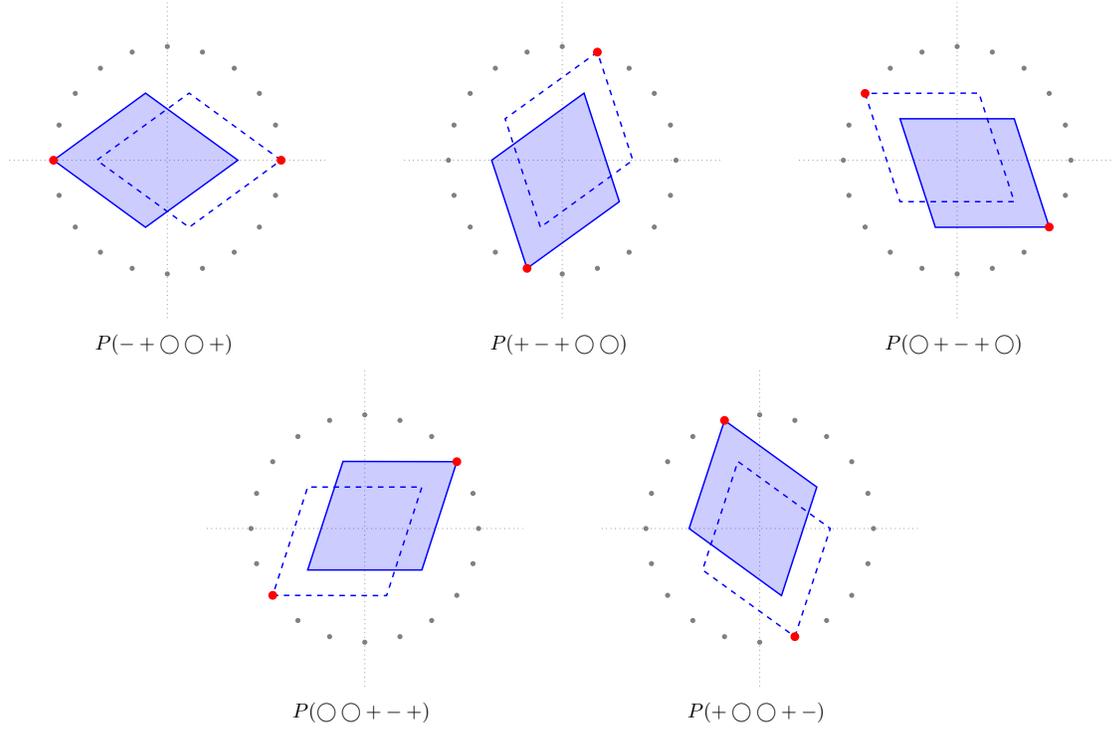

The construction works in general as follows. Whenever the cutting plane intersects a $k$-boundary of the Voronoi complex, the dual $(4{-}k)$-boundary is projected to the cutting plane. 
In our case, we wish to get the rhombic Penrose tiling. Therefore, we restrict ourselves to the skeletons $\mathscr{V}^{(2)}$ and $\mathscr{V}^{*(2)}$. Figure \ref{fig:projections_rhombi} shows the projections of the different (modulo translation) 2-boundaries into the $\mathcal{S}^{\parallel}$, which are the thick Penrose rhombuses. 

We choose as the cutting plane a translate of $\mathcal{S}^{\parallel}$ by a vector $\bs{c}_{{}_\perp} \in \mathcal{S}^{\perp}$. 
To ensure aperiodicity, we have to choose $\bs{c}_{{}_\perp}$ such that it is not contained in any $\pi_{{}_\perp}$-projection of any 1-boundary of $P^{*} \in \mathscr{V}^{*(2)}$. See \cite{BKSZ90} for further details.
The vector $\bs{c}_{{}_\perp}$ restricts the elements of $\mathscr{V}^{*(2)}$ which one projects on $\mathcal{S}^{\parallel}$, since the cutting plane  $\bs{c}_{{}_\perp} + \mathcal{S}^{\parallel}$ intersects 2-boundary $P$ if and only if $\pi_{{}_\perp}(P^{*})$ contains $\bs{c}_{{}_\perp}$. 
The resulting tiling (which depends on the choice of $\bs{c}_{{}_\perp}$) can be described as
\begin{equation}
\label{eq:tiling_descr}
\mathscr{T}_{\mathrm{Pen}}(\bs{c}_{{}_\perp}) \, = \, \left\{ \pi_{{}_\parallel}(P) \ : \ P \in \mathscr{V}^{(2)}, \ \bs{c}_{{}_\perp} \in \pi_{{}_\perp} (P^{*}) \right\}.  
\end{equation}

Vertices of $\mathscr{T}_{\mathrm{Pen}}(\bs{c}_{{}_\perp})$ are projections of vertex points of certain Voronoi domains $V(\bs{v})$ for some $\bs{v}\in A_4$. 
As already discussed above, these vertices are elements of $A^{*}_4 \backslash A_4 $ and are of four translation types, as characterised by the function $r$. The vertex points $\bs{v}^{*}$ split into four orbits with respect to the translation action of $A_4$. For each orbit, one can choose a~representative $\bs{v}^{*}_i$, for example 
\begin{equation} 
\label{eq:repre}
\bs{v}^{*}_1 \defeq \bs{a}_1,\qquad \bs{v}^{*}_2 \defeq \bs{a}_1+\bs{a}_3,\qquad  \bs{v}^{*}_3 \defeq -\bs{a}_1-\bs{a}_3, \qquad \bs{v}^{*}_4 \defeq -\bs{a}_1. 
\end{equation}
From the construction of $\mathscr{T}_{\mathrm{Pen}}(\bs{c}_{{}_\perp})$, we see that a point $\bs{v}^{*} \in A^{*}_4 \backslash A_4 $ is a preimage of a~vertex point in $\mathscr{T}_{\mathrm{Pen}}(\bs{c}_{{}_\perp})$ if and only if $\bs{v}^{*} \in P$ with $\bs{c}_{{}_\perp} \in \pi_{{}_\perp} (P^{*})$ for some  $P\in \mathscr{V}^{(2)}$. 
Note that if a point is an element of a $k$-boundary, the dual $(4{-}k)$-boundary lies in the dual cell of that point and vice versa. 
So, $\bs{v}^{*} \in P$ iff $P^{*} \subset V^{*}(\bs{v}^{*})$ with $V^{*}(\bs{v}^{*})$ being a~translate of a~dual 4D- cell of the form \eqref{eq:4D_dual_1}, or \eqref{eq:4D_dual_2}. 
Thus, $\pi_{{}_\parallel}(\bs{v}^{*})$ is a vertex in $\mathscr{T}_{\mathrm{Pen}}(\bs{c}_{{}_\perp})$ if and only if $\bs{c}_{{}_\perp} \in \pi_{{}_\perp}(V^{*}(\bs{v}^{*}))$.
Two points $\bs{v}^{*}_{1}$, $\bs{v}^{*}_{2}$ with $r(\bs{v}^{*}_{1}) =  r(\bs{v}^{*}_{2})$ can only differ by a~lattice vector. The choice of representatives \eqref{eq:repre} allows us to relate any point $\bs{v}^{*}$ with one of them. Define 
$q(\bs{v}^{*}) \defeq \bs{v}^{*} - \bs{v}^{*}_{r(\bs{v}^{*})} \in A_4$ for any $\bs{v}^{*}$. Since
\[ V^{*}(\bs{v}^{*}) = V^{*}(\bs{v}^{*} - \bs{v}^{*}_{r(\bs{v}^{*})} + \bs{v}^{*}_{r(\bs{v}^{*})}) = q(\bs{v}^{*}) + V^{*}(\bs{v}^{*}_{r(\bs{v}^{*})}), \]
one has 
\[ \bs{c}_{{}_\perp} \in \pi_{{}_\perp}(V^{*}(\bs{v}^{*})) \ \Longleftrightarrow \ \bs{c}_{{}_\perp} - \pi_{{}_\perp}(q(\bs{v}^{*})) \in \pi_{{}_\perp}\bigl( V^{*}(\bs{v}^{*}_{r(\bs{v}^{*})}) \bigr).\]
This allows us to rewrite the set of vertex points of $\mathscr{T}_{\mathrm{Pen}}(\bs{c}_{{}_\perp})$ as 
\begin{equation}
    \left\{ \pi_{{}_\parallel}(\bs{v}^{*}) \ : \ \bs{v}^{*} \in A^{*}_4 \backslash A_4, \ \bs{c}_{{}_\perp} - \pi_{{}_\perp}(q(\bs{v}^{*})) \in \pi_{{}_\perp}\bigl( V^{*}(\bs{v}^{*}_{r(\bs{v}^{*})}) \bigr)  \right\}.
\end{equation}

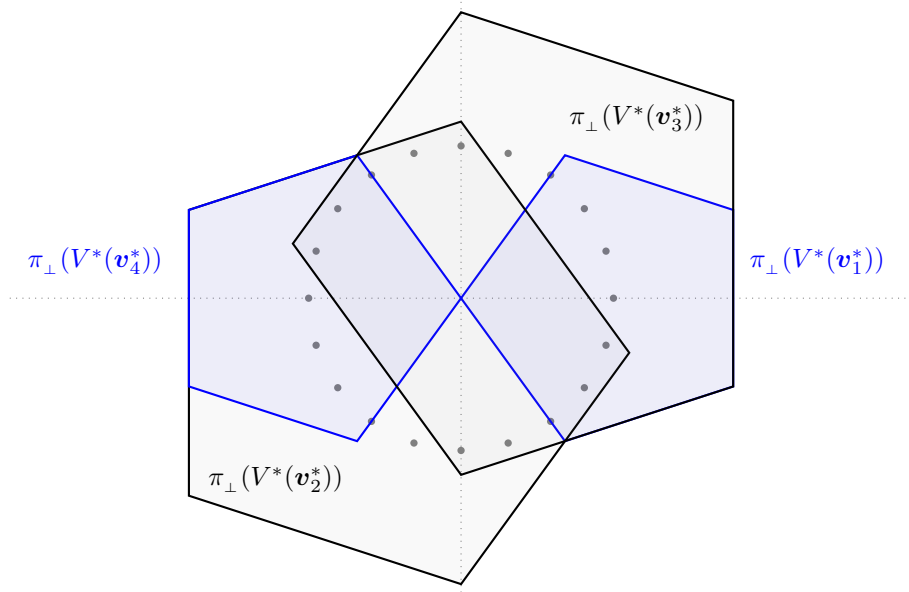
\begin{figure}
    \centering
    \scalebox{1}{
\begin{tikzpicture}
    \draw[thin, gray,dotted] (-6,0)--(6,0);
    \draw[thin, gray,dotted] (0,-4)--(0,4);
    \node[regular polygon, regular polygon sides=20, shape border rotate = -81, inner sep=40] (s) at (0,0) {} ;
    \foreach \x in {1,...,20}
			    \node[circle,draw=gray, fill=gray, inner sep=0pt,minimum size=2pt] (s.\x) at (s.corner \x) {};
    \draw[thick,blue, fill=blue, fill opacity = 0.05] (0,0) -- (1.3818,1.902)-- (3.618,1.174)-- (3.618,-1.174) -- (1.3818,-1.902) -- cycle;
    \draw[thick,black, fill=gray, fill opacity = 0.05] (0,2.35) -- (2.236,-0.726)-- (0,-3.804)-- (-3.618,-2.628) -- (-3.618,1.175) -- cycle;
    \draw[thick,blue, fill=blue, fill opacity = 0.05] (0,0) -- (-1.3818,-1.902)-- (-3.618,-1.174)-- (-3.618,1.174) -- (-1.3818,1.902) -- cycle;
    \draw[thick,black, fill=gray, fill opacity = 0.05] (0,-2.35) -- (-2.236,0.726)-- (0,3.804)-- (3.618,2.628) -- (3.618,-1.175) -- cycle;
    \node[align=left] (A) at (4.8,0.5) {$\color{blue}\small\pi_{{}_\perp}(V^{*}(\bs{v}^{*}_1))$ };
    \node[align=left] (B) at (-2.4,-2.4) {$\small\pi_{{}_\perp}(V^{*}(\bs{v}^{*}_2))$ };
    \node[align=left] (D) at (-4.8,0.5) {$\color{blue}\small\pi_{{}_\perp}(V^{*}(\bs{v}^{*}_4))$ };
    \node[align=left] (C) at (2.4,2.4) {$\small\pi_{{}_\perp}(V^{*}(\bs{v}^{*}_3))$ };
\end{tikzpicture}
}
\caption{Projections $\pi_{{}_\perp}\bigl( V^{*}(\bs{v}^{*}_{i})\bigr)\subset \mathcal{S}^{\perp}$ corresponding to the windows. The blue pentagons carry the $\pi_{{}_\perp}$-projections of shallow holes, whereas the black ones comprise the projections of deep holes. Note that for every window, there exists its own lattice. Thus even though there is a non-trivial intersection of windows, the resulting points must differ, as one expects. The gray points are the $20^{\mathrm{th}}$ roots of unity scaled by the factor $\sqrt{\tfrac{2}{5}}$.  }
    \label{fig:windows}
\end{figure}

This description shows that the set of vertices can be understood as 4 cut-and-project sets with lattices $\bs{v}^{*}_i + A_4 \subset \mathcal{S}$ and windows $\pi_{{}_\perp}\bigl( V^{*}(\bs{v}^{*}_{i})\bigr) \subset \mathcal{S}^{\perp}$, $1 \leqslant i \leqslant 4$. Figure \ref{fig:windows} shows all four windows $\pi_{{}_\perp}\bigl( V^{*}(\bs{v}^{*}_{i})\bigr)$ in $\mathcal{S}^{\perp}$, for more detail see Example 7.11. and Remark 7.8. in \cite{TAO}. 

Once we have established the description of all vertices of $\mathscr{T}_{\mathrm{Pen}}(\bs{c}_{{}_\perp})$, we can further determine a \emph{vertex configuration} of each vertex, i.e. all tiles in $\mathscr{T}_{\mathrm{Pen}}(\bs{c}_{{}_\perp})$ surrounding the vertex $\pi_{{}_\parallel}(\bs{v}^{*})$.
The description \eqref{eq:tiling_descr} provides us a characterization of the tiles surrounding $\pi_{{}_\parallel}(\bs{v}^{*})$. 
Indeed, a tile $\pi_{{}_\parallel}(P)$ belongs to a vertex configuration of $\pi_{{}_\parallel}(\bs{v}^{*})$ if and only if $P\in \mathscr{V}^{(2)}$, $\bs{v}^{*} \in P$ and $\bs{c}_{{}_\perp} \in \pi_{{}_\perp}(P^{*})$. 
The problem of finding a vertex configuration around an arbitrary vertex point can be reduced using translation symmetry. 
We can restrict ourselves to finding all vertex configurations around a representative of each translation class, i.e., around the points  $\bs{v}^{*}_i$. 
Then, we can rewrite the conditions above as $P^* -q(\bs{v}^{*}) \subset V^{*}(\bs{v}^*_{r(\bs{v}^*)})$ and $\bs{c}_{{}_\perp} - \pi_{{}_\perp}(q(\bs{v}^*)) \in \pi_{{}_\perp}(P^{*})- \pi_{{}_\perp}(q(\bs{v}^*))$. 
So, $P$ belongs to a vertex configuration of a point $\bs{v}^*$ if and only if it translate of its dual $P^{*}$ by $q(\bs{v}^*)$ is a 2-boundary of the dual cell $V^*(\bs{v}^*_{r(\bs{v}^*)})$. 
This gives an algorithm for obtaining the complete vertex configuration around the vertex $\pi_{{}_\parallel}(\bs{v}^{*})$ as follows.
\begin{enumerate}
    \item Find all $\bs{w}^* \in A_4^* \backslash A_4$ such that $\bs{c}_{{}_\perp} - \pi_{{}_\perp}(q(\bs{w}^*)) \in \pi_{{}_\perp}(V^{*}(\bs{v}^*_{r(\bs{v}^*)})$.
    \item For all $\bs{w}^*$ found in step 1, take the 2-boundary $P^{*}$ of the dual cell $V^{*}(\bs{v}^*_{r(\bs{v}^*)})$ with $\bs{c}_{{}_\perp} - \pi_{{}_\perp}(q(\bs{w}^*)) \in \pi_{{}_\perp}(P^{*})$. Then, $\pi_{{}_\parallel}(\bs{w}^{*}) + \pi_{{}_\parallel}(P)$ is a tile around $\pi_{{}_\parallel}(\bs{v}^{*})$.
\end{enumerate}

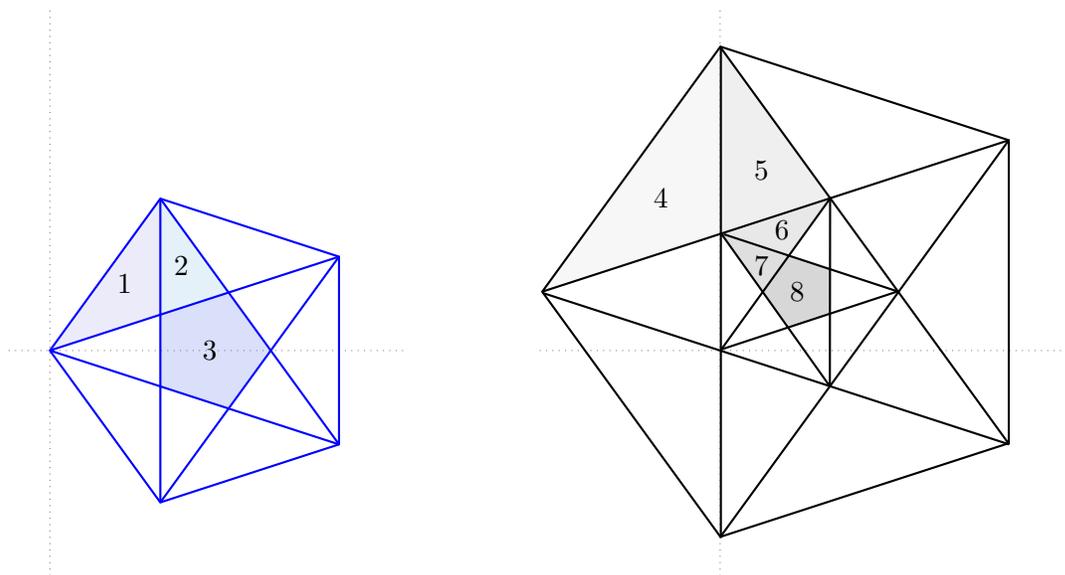
\begin{figure}
    \centering
    \begin{tikzpicture}[x=0.75pt,y=0.75pt,yscale=-1.3,xscale=1.3]

\draw[thin, gray,dotted]    (140,218) -- (294.15,218) ;
\draw[thin, gray,dotted]   (156,86) -- (156,306) ;

\draw  [draw opacity=0][fill={rgb, 255:red, 74; green, 100; blue, 226 }  ,fill opacity=0.2 ] (241.68,218) -- (234.97,227.23) -- (225.32,240.52) -- (198.84,231.92) -- (198.84,204.08) -- (225.32,195.48) -- (241.68,218) -- cycle ;
\draw  [draw opacity=0][fill={rgb, 255:red, 120; green, 188; blue, 220 }  ,fill opacity=0.2 ] (198.84,159.03) -- (225.32,195.48) -- (198.84,204.08) -- (198.84,159.03) -- cycle ;
\draw  [draw opacity=0][fill={rgb, 255:red, 155; green, 155; blue, 230}  ,fill opacity=0.2 ] (156,218) -- (156,218) -- (198.84,159.03) -- (198.84,204.08) -- (156,218) -- cycle ;

\draw[thick, blue]  (156,218) -- (198.84,159.03) -- (268.16,181.56) -- (268.16,254.44) -- (198.84,276.97) -- cycle ;
\draw[thick,blue]    (198.84,159.03) -- (268.16,254.44) ;
\draw[thick,blue]   (198.84,276.97) -- (198.84,159.03) ;
\draw[thick,blue]    (156,218) -- (268.16,254.44) ;
\draw[thick,blue]    (156,218) -- (268.16,181.56) ;
\draw[thick,blue]    (198.84,276.97) -- (268.16,181.56) ;

\node (1) at (185,192) {$1$};
\node (2) at (207,185) {$2$};
\node (3) at (218,218) {$3$};

\draw[thin, gray,dotted]    (416,86) -- (416,306)  ;
\draw[thin, gray,dotted]   (346,218) -- (554,218) ;
\draw  [draw opacity=0][fill={rgb, 255:red, 59; green, 59; blue, 59 }  ,fill opacity=0.2 ] (432.51,195.1) -- (442.51,181.33) -- (458.7,186.59) -- (458.7,203.61) -- (442.51,208.87) -- cycle ;
\draw  [draw opacity=0][fill={rgb, 255:red, 96; green, 96; blue, 96 }  ,fill opacity=0.2 ] (442.69,180.92) -- (432.38,195.08) -- (416.47,172.63) -- cycle ;
\draw  [draw opacity=0][fill={rgb, 255:red, 140; green, 140; blue, 140 }  ,fill opacity=0.2 ] (458.68,158.95) -- (442.69,180.92) -- (416.47,172.63) -- cycle ;
\draw  [draw opacity=0][fill={rgb, 255:red, 178; green, 178; blue, 178 }  ,fill opacity=0.2 ] (416.29,100.46) -- (458.68,158.95) -- (416.47,172.63) -- cycle ;
\draw  [draw opacity=0][fill={rgb, 255:red, 215; green, 215; blue, 215 }  ,fill opacity=0.2 ] (416.29,100.46) -- (416.47,172.63) -- (346.99,195.16) -- cycle ;

\draw[thick,black]   (346.99,195.16) -- (416.13,100) -- (528,136.35) -- (528,253.98) -- (416.13,290.33) -- cycle ;
\draw[thick,black]     (416.29,100.46) -- (528.16,254.44) ;
\draw[thick,black]     (416.29,100.46) -- (416.29,290.79) ;
\draw[thick,black]     (528,136.35) -- (346.99,195.16) ;
\draw[thick,black]     (528,136.35) -- (416.13,290.33) ;
\draw[thick,black]     (347.15,195.63) -- (528.16,254.44) ;
\draw[thick,black]     (416.47,172.63) -- (458.58,231.89) ;
\draw[thick,black]     (458.68,158.95) -- (458.58,231.89) ;
\draw[thick,black]     (458.68,158.95) -- (416.47,217.58) ;
\draw[thick,black]     (484.89,195.05) -- (417.42,172.84) ;
\draw[thick,black]     (484.89,195.05) -- (416.47,217.58) ;

\node (4) at (393,159) {$4$};
\node (5) at (432,148) {$5$};
\node (6) at (440,171) {$6$};
\node (7) at (432,185) {$7$};
\node (8) at (446,195) {$8$};

\end{tikzpicture}  
    \caption{Subdivision of $\pi_{{}_\perp}(V^{*}(\bs{v}^*_{1}))$ (blue) and $\pi_{{}_\perp}(V^{*}(\bs{v}^*_{3}))$ (black) into elementary polygons. The eight possible vertex configurations (modulo rotation by $ \tfrac{2\pi}{5}$ and space inversion) correspond to eight distinct elementary polygons.  }
    \label{fig:elem_poly}
\end{figure}

\begin{figure}
    \centering
     

\begin{tikzpicture}[x=0.75pt,y=0.75pt,yscale=-1,xscale=1]
%1
\draw  [color=black  ,draw opacity=1 ][fill=blue  ,fill opacity=0.3 ] (94.53,79.46) -- (86.48,104.18) -- (78.46,79.45) -- (86.51,54.73) -- cycle ;
\draw  [color=black  ,draw opacity=1 ][fill=blue  ,fill opacity=0.3 ] (78.46,79.45) -- (57.42,94.72) -- (65.47,70) -- (86.51,54.73) -- cycle ;
\draw  [color=black  ,draw opacity=1 ][fill=red  ,fill opacity=0.3 ] (65.44,119.46) -- (57.42,94.72) -- (78.46,79.45) -- (86.48,104.18) -- cycle ;
\node[circle,draw=black, fill=black, inner sep=1pt, minimum size= 4pt] (1) at (78.46,80) {};

%2
\draw  [color=black  ,draw opacity=1 ][fill=red  ,fill opacity=0.3 ] (234.05,55.28) -- (226,80) -- (200,79.99) -- (208.05,55.26) -- cycle ;
\draw  [color=black  ,draw opacity=1 ][fill=red  ,fill opacity=0.3 ] (186.99,120) -- (178.97,95.27) -- (200.01,80) -- (208.03,104.73) -- cycle ;
\draw  [color=black  ,draw opacity=1 ][fill=red  ,fill opacity=0.3 ] (234.03,104.74) -- (208.03,104.73) -- (200.01,80) -- (226.01,80.01) -- cycle ;
\draw  [color=black  ,draw opacity=1 ][fill=blue  ,fill opacity=0.3 ] (200,79.99) -- (178.96,95.26) -- (187.01,70.53) -- (208.05,55.26) -- cycle ;
\node[circle,draw=black, fill=black, inner sep=1pt, minimum size= 4pt] (2) at (200,80) {};

%3
\draw  [color=black  ,draw opacity=1 ][fill=red  ,fill opacity=0.3 ] (354.05,55.28) -- (346,80) -- (320,79.99) -- (328.05,55.26) -- cycle ;
\draw  [color=black  ,draw opacity=1 ][fill=red  ,fill opacity=0.3 ] (320,79.99) -- (298.96,95.26) -- (277.93,79.96) -- (298.97,64.69) -- cycle ;
\draw  [color=black  ,draw opacity=1 ][fill=red  ,fill opacity=0.3 ] (320,79.99) -- (298.97,64.69) -- (307.02,39.97) -- (328.05,55.26) -- cycle ;
\draw  [color=black  ,draw opacity=1 ][fill=red  ,fill opacity=0.3 ] (306.98,119.99) -- (298.96,95.26) -- (320,79.99) -- (328.02,104.72) -- cycle ;
\draw  [color=black  ,draw opacity=1 ][fill=red  ,fill opacity=0.3 ] (354.02,104.73) -- (328.02,104.72) -- (320,79.99) -- (346,80) -- cycle ;
\node[circle,draw=black, fill=black, inner sep=1pt, minimum size= 4pt] (3) at (320,80) {};

%4
\draw  [color=black  ,draw opacity=1 ][fill=red  ,fill opacity=0.3 ] (448,55.28) -- (439.95,80) -- (413.95,79.99) -- (422,55.26) -- cycle ;
\draw  [color=black  ,draw opacity=1 ][fill=red  ,fill opacity=0.3 ] (460.98,95.29) -- (439.95,80) -- (448,55.28) -- (469.03,70.57) -- cycle ;
\draw  [color=black  ,draw opacity=1 ][fill=blue  ,fill opacity=0.3 ] (434.98,95.28) -- (413.95,79.99) -- (439.95,80) -- (460.98,95.29) -- cycle ;
\node[circle,draw=black, fill=black, inner sep=1pt, minimum size= 4pt] (4A) at (448,55.28) {};
\node[circle,draw=black, fill=black, inner sep=1pt, minimum size= 4pt] (4B) at (434.98,95.28) {};

%5
\draw  [color=black  ,draw opacity=1 ][fill=blue  ,fill opacity=0.3 ] (88.03,175.28) -- (79.99,200) -- (71.97,175.27) -- (80.01,150.55) -- cycle ;
\draw  [color=black  ,draw opacity=1 ][fill=red  ,fill opacity=0.3 ] (101.01,215.3) -- (79.97,230.57) -- (58.94,215.28) -- (79.99,200) -- cycle ;
\draw  [color=black  ,draw opacity=1 ][fill=red  ,fill opacity=0.3 ] (101.01,215.3) -- (79.99,200) -- (88.03,175.28) -- (109.06,190.57) -- cycle ;
\draw  [color=black  ,draw opacity=1 ][fill=blue  ,fill opacity=0.3 ] (58.94,215.28) -- (32.95,215.26) -- (53.99,199.99) -- (79.99,200) -- cycle ;
\draw  [color=black  ,draw opacity=1 ][fill=red  ,fill opacity=0.3 ] (79.99,200) -- (53.99,199.99) -- (45.97,175.26) -- (71.97,175.27) -- cycle ;
\node[circle,draw=black, fill=black, inner sep=1pt, minimum size= 4pt] (5A) at (45.97,175.26) {};
\node[circle,draw=black, fill=black, inner sep=1pt, minimum size= 4pt] (5B) at (88.03,175.28) {};
\node[circle,draw=black, fill=black, inner sep=1pt, minimum size= 4pt] (5C) at (58.94,215.28) {};

%6
\draw  [color=black  ,draw opacity=1 ][fill=blue  ,fill opacity=0.3 ] (208.03,175.28) -- (199.99,200) -- (191.97,175.27) -- (200.01,150.55) -- cycle ;
\draw  [color=black  ,draw opacity=1 ][fill=red  ,fill opacity=0.3 ] (242.06,200.03) -- (221.01,215.3) -- (199.99,200) -- (221.03,184.73) -- cycle ;
\draw  [color=black  ,draw opacity=1 ][fill=red  ,fill opacity=0.3 ] (212.97,240.02) -- (191.94,224.73) -- (199.99,200) -- (221.01,215.3) -- cycle ;
\draw  [color=black  ,draw opacity=1 ][fill=blue  ,fill opacity=0.3 ] (178.94,215.28) -- (152.95,215.26) -- (173.99,199.99) -- (199.99,200) -- cycle ;
\draw  [color=black  ,draw opacity=1 ][fill=red  ,fill opacity=0.3 ] (199.99,200) -- (173.99,199.99) -- (165.97,175.26) -- (191.97,175.27) -- cycle ;
\draw  [color=black  ,draw opacity=1 ][fill=blue  ,fill opacity=0.3 ] (221.03,184.73) -- (199.99,200) -- (208.03,175.28) -- (229.08,160.01) -- cycle ;
\draw  [color=black  ,draw opacity=1 ][fill=blue  ,fill opacity=0.3 ] (191.94,224.73) -- (170.9,240) -- (178.94,215.28) -- (199.99,200) -- cycle ;
\node[circle,draw=black, fill=black, inner sep=1pt, minimum size= 4pt] (6A) at (165.97,175.26) {};
\node[circle,draw=black, fill=black, inner sep=1pt, minimum size= 4pt] (6B) at (208.03,175.28) {};
\node[circle,draw=black, fill=black, inner sep=1pt, minimum size= 4pt] (6C) at (178.94,215.28)  {};
\node[circle,draw=black, fill=black, inner sep=1pt, minimum size= 4pt] (6D) at (242.06,200.03)  {};
\node[circle,draw=black, fill=black, inner sep=1pt, minimum size= 4pt] (6E) at (212.97,240.02) {};

%7
\draw  [color=black  ,draw opacity=1 ][fill=red  ,fill opacity=0.3 ] (319.92,199.99) -- (311.88,224.71) -- (285.88,224.7) -- (293.92,199.98) -- cycle ;
\draw  [color=black  ,draw opacity=1 ][fill=blue  ,fill opacity=0.3 ] (327.97,175.27) -- (319.92,199.99) -- (311.9,175.26) -- (319.95,150.53) -- cycle ;
\draw  [color=black  ,draw opacity=1 ][fill=red  ,fill opacity=0.3 ] (361.99,200.01) -- (340.95,215.28) -- (319.92,199.99) -- (340.97,184.72) -- cycle ;
\draw  [color=black  ,draw opacity=1 ][fill=red  ,fill opacity=0.3 ] (332.9,240.01) -- (311.88,224.71) -- (319.92,199.99) -- (340.95,215.28) -- cycle ;
\draw  [color=black  ,draw opacity=1 ][fill=red  ,fill opacity=0.3 ] (319.92,199.99) -- (293.92,199.98) -- (285.9,175.24) -- (311.9,175.26) -- cycle ;
\draw  [color=black  ,draw opacity=1 ][fill=blue  ,fill opacity=0.3 ] (340.97,184.72) -- (319.92,199.99) -- (327.97,175.27) -- (349.01,159.99) -- cycle ;
\node[circle,draw=black, fill=black, inner sep=1pt, minimum size= 4pt] (7A) at (285.97,175.26) {};
\node[circle,draw=black, fill=black, inner sep=1pt, minimum size= 4pt] (7B) at (328.03,175.28) {};
\node[circle,draw=black, fill=black, inner sep=1pt, minimum size= 4pt] (7C) at (285.88,224.7)  {};
\node[circle,draw=black, fill=black, inner sep=1pt, minimum size= 4pt] (7D) at (362.06,200.03)  {};
\node[circle,draw=black, fill=black, inner sep=1pt, minimum size= 4pt] (7E) at (332.97,240.02) {};

%8
\draw  [color=black  ,draw opacity=1 ][fill=red  ,fill opacity=0.3 ] (439.92,199.99) -- (431.9,175.26) -- (452.94,159.99) -- (460.97,184.72) -- cycle ;
\draw  [color=black  ,draw opacity=1 ][fill=red  ,fill opacity=0.3 ] (439.92,199.99) -- (431.88,224.71) -- (405.88,224.7) -- (413.92,199.98) -- cycle ;
\draw  [color=black  ,draw opacity=1 ][fill=red  ,fill opacity=0.3 ] (481.99,200.01) -- (460.95,215.28) -- (439.92,199.99) -- (460.97,184.72) -- cycle ;
\draw  [color=black  ,draw opacity=1 ][fill=red  ,fill opacity=0.3 ] (452.9,240.01) -- (431.88,224.71) -- (439.92,199.99) -- (460.95,215.28) -- cycle ;
\draw  [color=black  ,draw opacity=1 ][fill=red  ,fill opacity=0.3 ] (439.92,199.99) -- (413.92,199.98) -- (405.9,175.24) -- (431.9,175.26) -- cycle ;
\node[circle,draw=black, fill=black, inner sep=1pt, minimum size= 4pt] (7A) at (405.97,175.26) {};
\node[circle,draw=black, fill=black, inner sep=1pt, minimum size= 4pt] (7B) at (452.94,159.99) {};
\node[circle,draw=black, fill=black, inner sep=1pt, minimum size= 4pt] (7C) at (405.88,224.7)  {};
\node[circle,draw=black, fill=black, inner sep=1pt, minimum size= 4pt] (7D) at (482.06,200.03)  {};
\node[circle,draw=black, fill=black, inner sep=1pt, minimum size= 4pt] (7E) at (452.97,240.02) {};

\node[align=center] (p1) at (45,55) {$1$ };
\node[align=center] (p2) at (165,55) {$2$ };
\node[align=center] (p3) at (285,55) {$3$ };
\node[align=center] (p4) at (405,55) {$4$ };
\node[align=center] (p5) at (45,155) {$5$ };
\node[align=center] (p6) at (165,155) {$6$ };
\node[align=center] (p7) at (285,155) {$7$ };
\node[align=center] (p8) at (405,155) {$8$ };

\end{tikzpicture}
    \caption{All possible vertex configurations (up to rotation by $ \tfrac{2\pi}{5}$ and space inversion)  which are in one-to-one correspondence with the elementary polygons in Figure~\ref{fig:elem_poly}. The black points indicate the positions of shallow holes. }
    \label{fig:vertex_conf}
\end{figure}
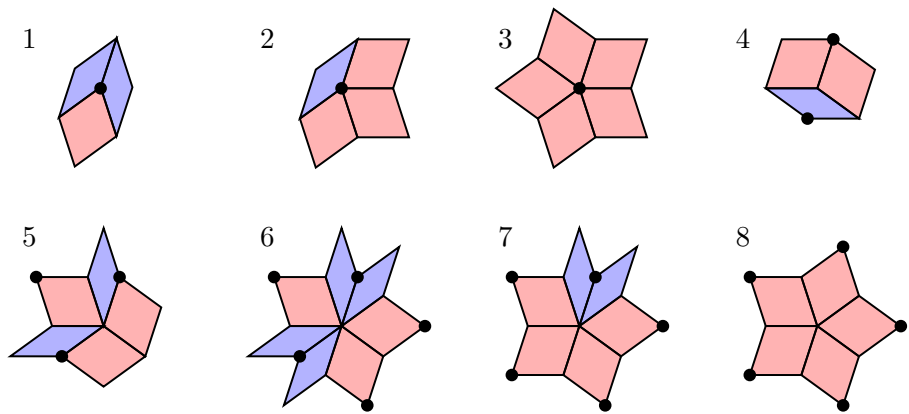

We chose $\bs{c}_{{}_\perp}$ so that $\bs{c}_{{}_\perp} - \pi_{{}_\perp}(q(\bs{w}^*))$  lies in the interior of $ \pi_{{}_\perp}(P^{*})$. 
This is a crucial observation.
It forces all tiles $\pi_{{}_\parallel}(\bs{w}^{*}) + \pi_{{}_\parallel}(P)$ belonging to a particular vertex configuration to have, at the level of $\pi_{{}_\perp}(P^{*})$, an overlap in the $\pi_{{}_\perp}(V^{*}(\bs{v}^*_{r(\bs{v}^*)})$. 
We can use this property to determine and characterise all possible vertex configurations with respect to translations in $\mathcal{S}^{\perp}$ as follows: 
A set $\mathscr{W}\subset \mathscr{V}^{(2)}$ of 2-boundaries is a~valid vertex configuration of a~vertex of type $i$ if and only if $\mathscr{W}$ is maximal with respect to the property that $\bigcap_{P\in \mathscr{W}} \pi_{{}_\perp}(P^*)$ is non-empty. 
The projection of 2-boundaries of the dual cells $V^{*}(\bs{v}^*_{i})$ divides the $\pi_{{}_\perp}(V^{*}(\bs{v}^*_{i}))$ into convex polygons, so-called \emph{elementary polygons} \cite{BKSZ90}. They have pairwise distinct interiors, each representing a distinct vertex configuration (and vice versa). Figure \ref{fig:elem_poly} shows the elementary polygons for $\pi_{{}_\perp}(V^{*}(\bs{v}^*_{1}))$ and $\pi_{{}_\perp}(V^{*}(\bs{v}^*_{3}))$. The corresponding vertex configurations are shown in Figure \ref{fig:vertex_conf}.

The choice of the cutting plane ensures that the projections of the vertices of a valid infinite Penrose tiling into $\mathcal{S}^{\perp}$ are dense and uniformly distributed. 
Thus, we can use them to determine the frequencies of the vertex configurations via the areas of elementary polygons. Denote by $E$ an elementary polygon. 
Then, the relative frequency $\nu_{\mathscr{C}(E)}$ of vertex configuration $\mathscr{C}(E)$ corresponding to the elementary polygon $E$ is given by  
\begin{equation}
    \label{eq:freq_element}
    \nu_{\mathscr{C}(E)} \, = \, \myfrac{\mathrm{Area}(E)}{\sum_{i=1}^4\mathrm{Area}(\pi_{{}_\perp}(V^{*}(\bs{v}^*_{i}))) } = \myfrac{\mathrm{Area}(E)}{\sqrt{5\bigl(5+2\sqrt{5}\bigr)}},
\end{equation}
i.e., exactly by the fraction of the total area of windows it occupies. We list the frequencies of all vertex configurations. We include the frequency of given patch as well as the cumulative frequency of all patches of the same type, i.e., all patches that lie in the same orbit under the rotation and space inversion. 

\begin{table}[h!]
\centering
\bgroup
\def\arraystretch{1.4}
\begin{tabular}{c||c|c}
     Vertex config. & Frequency $\nu_{\mathscr{C}(E_i)}$  & Total frequency $\nu_i$ \\
     \hline \hline
     1 & $\tfrac{1}{10}(5-3\tau)$ & $5-3\tau = \tau^{-4}$ \\
     2 & $\tfrac{1}{10}(5\tau-8)$ & $5\tau-8= \tau^{-5}$ \\
     3 & $\tfrac{1}{10}(18-11\tau)$ & $\tfrac{1}{5}(18-11\tau)= \tfrac{2\tau-1}{5}\tau^{-5}$ \\
     4 & $\tfrac{1}{10}(2-\tau)$ & $2-\tau= \tau^{-2}$ \\
     5 & $\tfrac{1}{10}(2\tau-3)$ & $2\tau-3= \tau^{-3}$ \\
     6 & $\tfrac{1}{10}(13 - 8\tau)$ & $13-8\tau= \tau^{-6}$ \\
     7 & $\tfrac{1}{10}(13\tau-21)$ & $13\tau-21= \tau^{-7}$ \\
     8 & $\tfrac{1}{10}(47 - 29\tau)$ & $\tfrac{1}{5}(47 - 29\tau)= \tfrac{2\tau-1}{5}\tau^{-7}$ \\
     \hline \hline
\end{tabular}
\vspace{0.5cm}
\egroup
\caption{Frequencies of vertex configurations in Penrose tilings, all belonging to $\QQ(\tau)$ with $\tau$ being the golden ratio. The second column shows frequencies of particular patches, those in Figure \ref{fig:vertex_conf}. The last column gives the total frequencies of a patch of a given type, i.e., a patch and all its images under the allowed rotations and space inversion. }
\end{table}

The sum of all total frequencies equals  one; thus, we get a consistency check. Since Penrose tiling defines a strictly ergodic dynamical system (in the usual way) \cite{R96}, we conclude that there are no other vertex configurations. If there were any others, they would come with a~strictly positive measure, which is the patch frequency. 

Recall that the \emph{frequency module} of a tiling space $\XX$ (in our case, the tiling space generated by rhombic Penrose tilings) is the minimal $\ZZ$-module $\mathcal{M}_{\XX}$ that contains all frequencies of finite patches of the tiling. Here, we get the following specific result. 

\begin{prop}
The frequency module $\mathcal{M}_{\mathscr{T}_{\mathrm{Pen}}}$ of the Penrose tiling is 
\begin{equation}
    \mathcal{M}_{\mathscr{T}_{\mathrm{Pen}}} \, = \, \myfrac{1}{10}\ZZ[\tau].
\end{equation}
\end{prop}
\noindent \textit{Proof.} Consider any patch of the Penrose tiling. 
We can always find an $n\in \NN$ such that this given patch is contained in a level-$n$ supertile of some vertex configuration $\mathscr{C}$. 
Since the Penrose tiling is an inflation/deflation tiling, its level-$n$ supertiles around given vertex configuration are equivalent to the original vertex configuration scaled by factor $\tau^n$. 
Therefore, the supertile itself has a frequency given by $\tfrac{1}{\tau^{2n}}\cdot \nu_{\mathscr{C}}$. The factor $\tfrac{1}{\tau^{2n}}$ comes from the observation that the frequency is inversely proportional to the area.
Since $\tau$ is a unit in $\ZZ[\tau]$, so is $\tau^{2n}$. 
Thus, to determine the frequency module, it suffices to consider the $\ZZ$-module generated by $\nu_{\mathscr{C}(E)_i}$, $1\leqslant i \leqslant 8$, i.e., 
\[ \mathcal{M}_{\mathscr{T}_{\mathrm{Pen}}}  = \langle \nu_{\mathscr{C}(E_i)}  \ : \ 1\leqslant i \leqslant 8 \rangle_{\ZZ}.  \]
Since $\nu_{\mathscr{C}(E_1)}+ 2\nu_{\mathscr{C}(E_5)}+\nu_{\mathscr{C}(E_4)} = \tfrac{1}{10}$ and $\nu_{\mathscr{C}(E_4)} + \nu_{\mathscr{C}(E_5)} = \tfrac{\tau -1}{10}$, one has 
\pushQED{\qed}
\[ \mathcal{M}_{\mathscr{T}_{\mathrm{Pen}}}  = \langle \nu_{\mathscr{C}(E_i)}  \ : \ 1\leqslant i \leqslant 8 \rangle_{\ZZ} = \myfrac{1}{10}\langle 1,\tau \rangle_{\ZZ} = \myfrac{1}{10} \ZZ[\tau] . \qedhere \]
\popQED

\section{General patch frequencies and their calculation}
\label{sec:freq}

The idea behind the above construction can be extended to any patch in Penrose rhombic tilings.
Choose a vertex of a tile and relate all tiles of the patch to this vertex. 
One has to be careful and consistently distinguish between deep and shallow holes. Then, one obtains a list of all tiles and their relative positions with respect to the chosen central tile. 
By transitioning into $\mathcal{S}^{\perp}$, one gets a list of all dual triangles and their relative distance. 
Their intersection determines the frequency of the patch in the same way as in the case of vertex configuration. 
This intersection is always a convex polygon (since one intersects a finite number of triangles), and its area can be computed easily. Note that some minimal subset of the triangles entirely determines this intersection, and working only with them can increase the computational speed considerably. 

\begin{figure}
    \centering
    \input{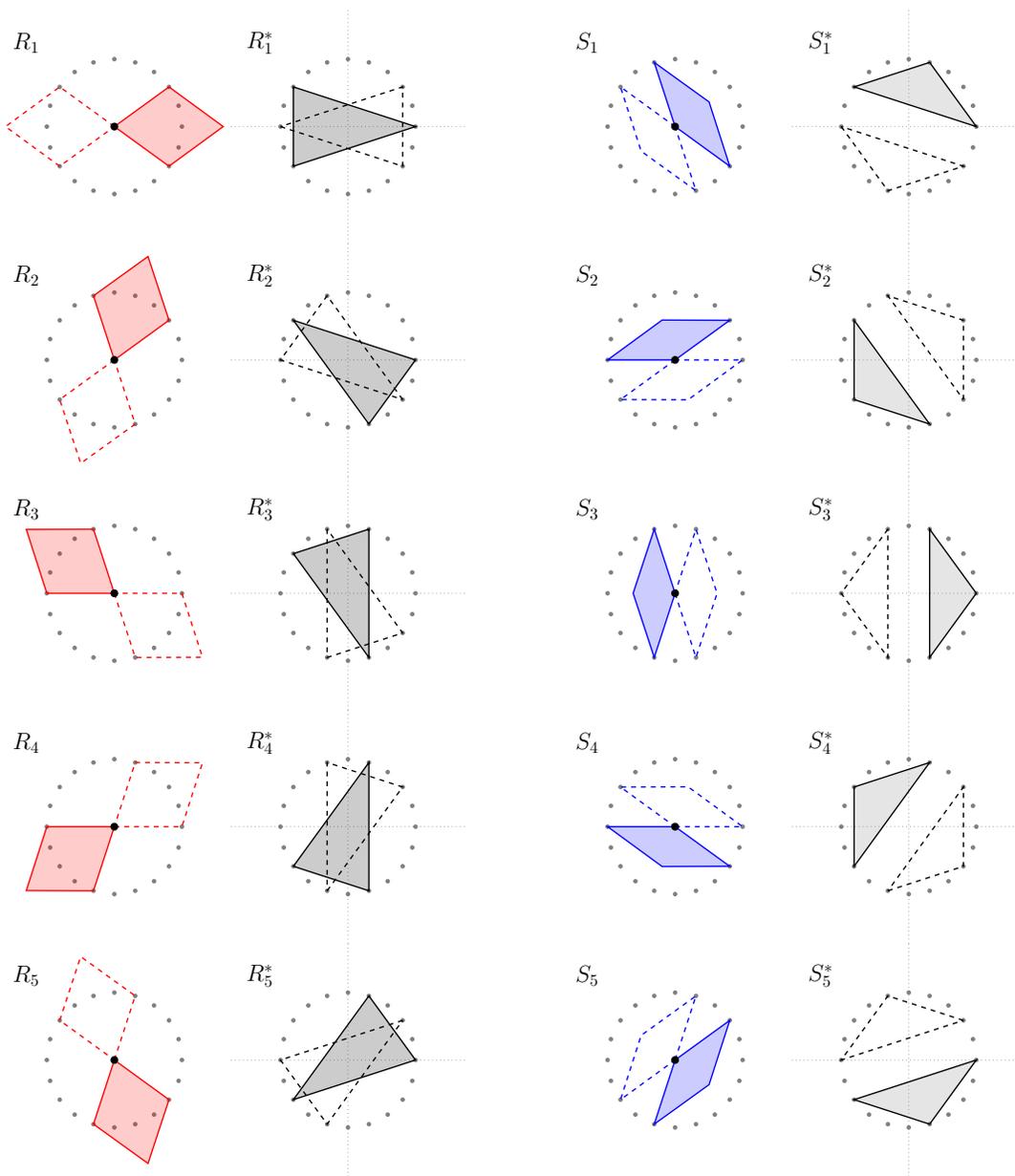}
    \caption{List of all possible tiles (with respect to their orientations and placement of a shallow hole) in rhombic Penrose tilings and their duals in $\mathcal{S}^{\perp}$. Tiles are depicted relative to the shallow hole. The exact correspondence between a tile in the list and a~projection of a 2-boundary is, for example, the following. If a tile of type $R_1$ corresponds to $\pi_{{}_\parallel}(P(-+\bigcirc \bigcirc +))$, the dual triangle $R^{*}_1$ \emph{is equal} to $\pi_{{}_\perp}(P^{*}(-+\bigcirc \bigcirc +)) + \pi_{{}_\perp}(\bs{a}_1)$, i.e., we capture its actual position in $\perp$-space. Fixing the positions of the duals allows us to work in coordinates relative to a given point, the ``origin". Then, everything is shifted by a~suitable vector representing the relative distances of objects to the ``origin". }
    \label{fig:vertex_list}
\end{figure}

Let us list, in Figure \ref{fig:vertex_list}, all possible tiles together with the shallow holes attached to each of them. We place them so that the shallow hole indicates the `origin' relative to the given tile. More precisely, we depict them in coordinates which are translated by the shallow vertex of a given tile. We also include the dual tile and its projection in $\mathcal{S}^{\perp}$. The projection is also centred on the relative origin. There is an extra advantage of such a choice, namely, the vertices of dual triangles in $\mathcal{S}^{\perp}$ are placed at the $20^{\mathrm{th}}$ roots of unity scaled by the factor $\sqrt{\tfrac{2}{5}}$. Since the frequency is given by a ratio of two areas, the scaling factor does not play a role. This allows a~precise calculation, simply by employing a suitable subfield of $\QQ(\ee^{\frac{\pi \ii}{10}})$. In fact, one can work with integer coefficients. 

We can now describe the algorithm that allows us to determine the frequency of a given patch. We start with an arbitrary finite patch of the Penrose tiling. 

\begin{enumerate}
    \item Detect all shallow holes in the patch. This can be done via the allowed vertex configurations. 
    \item Identify the type of each tile in the patch as $R_i$, $S_i$ or their space inversions ($RI_i$, $SI_i$) from list \ref{fig:vertex_list}.
    \item Choose any shallow hole, the ``origin", from the vertices of the patch and fix it. 
    \item Make a list of all positions of all tiles (their shallow holes) relative to the ``origin". Since the edges of the rhombuses are projections $\pi_{{}_\parallel}(\bs{a}_i)$, the resulting position vector can always be written as $\sum_{i\in \mathcal{I}}\epsilon_i \pi_{{}_\parallel}(\bs{a}_i)$ where $\mathcal{I}$ parametrises the path on the edges from the ``origin" to the desired point $\epsilon_i \in \{ \pm 1 \}$ denotes the orientation of the vectors $\pi_{{}_\parallel}(\bs{a}_i)$ in the path. 
    \item Apply the dual correspondence, i.e., to each translated tile $T + \sum_{i\in \mathcal{I}}\epsilon_i \pi_{{}_\parallel}(\bs{a}_i)$ assign the dual $T^{*} + \sum_{i\in \mathcal{I}}\epsilon_i \pi_{{}_\perp}(\bs{a}_i)$, with $T \in \{R_i, S_i, RI_i, SI_i \}$.
    \item Find an intersection of all  $T^{*} + \sum_{i\in \mathcal{I}}\epsilon_i \pi_{{}_\perp}(\bs{a}_i)$ from the list. This can be done via any clipping algorithm, for example the Sutherland–Hodgman algorithm \cite{SH74}. 
    \item Calculate the area of the intersection.
    \item Divide the area of the intersection by the total area of the windows, i.e., with $\sqrt{5\bigl(5+2\sqrt{5}\bigr)}$. This yields the relative frequency.
\end{enumerate}
Note that one can choose any clipping algorithm since one has to deal with triangles only (for different lattices one obtains general convex polygons).  Under this condition, most clipping algorithms are sufficiently robust. 
Moreover, at least the Sutherland–Hodgman algorithm ensures that the resulting coordinates of vertices of the intersection are contained in the same field as the coordinates of the polygons since each step of the algorithm relies on solving systems of two linear equations with coefficients being the coordinates of the vertices of the polygons.

Finally, computing the area of a polygon determined by its vertices can be done via the shoelace formula (or Gauss's area formula) \cite[p. 125]{KoecherKrieg}, which is also within the field. 

\begin{figure}[h!]
 \begin{subfigure}{.5\textwidth}
  \centering
  \scalebox{1.3}{
\begin{tikzpicture}[x=0.75pt,y=0.75pt,yscale=-1,xscale=1]

\draw  [color=black ,draw opacity=1 ][fill=red  ,fill opacity=0.35 ] (388.05,116.78) -- (380,141.5) -- (354,141.49) -- (362.05,116.76) -- cycle ;
\draw  [color=black ,draw opacity=1 ][fill=red  ,fill opacity=0.35 ] (354,141.49) -- (332.97,126.19) -- (341.02,101.47) -- (362.05,116.76) -- cycle ;
\draw  [color=black ,draw opacity=1 ][fill=red  ,fill opacity=0.35 ] (340.98,181.49) -- (332.96,156.76) -- (354,141.49) -- (362.02,166.22) -- cycle ;
\draw  [color=black ,draw opacity=1 ][fill=red  ,fill opacity=0.35 ] (388.02,166.23) -- (362.02,166.22) -- (354,141.49) -- (380,141.5) -- cycle ;
\draw  [color=black ,draw opacity=1 ][fill=red  ,fill opacity=0.35 ] (293.98,116.73) -- (285.96,92) -- (307,76.73) -- (315.02,101.46) -- cycle ;
\draw  [color=black ,draw opacity=1 ][fill=red  ,fill opacity=0.35 ] (264.89,156.72) -- (256.84,181.45) -- (230.84,181.43) -- (238.89,156.71) -- cycle ;
\draw  [color=black ,draw opacity=1 ][fill=red  ,fill opacity=0.35 ] (264.91,126.16) -- (238.91,126.14) -- (230.88,101.41) -- (256.88,101.43) -- cycle ;
\draw  [color=black ,draw opacity=1 ][fill=red  ,fill opacity=0.35 ] (306.93,206.2) -- (285.91,190.91) -- (293.95,166.18) -- (314.98,181.48) -- cycle ;

\draw  [color=black ,draw opacity=1 ][fill=blue  ,fill opacity=0.35 ] (272.91,181.45) -- (264.86,206.18) -- (256.84,181.45) -- (264.89,156.72) -- cycle ;
\draw  [color=black ,draw opacity=1 ][fill=blue  ,fill opacity=0.35 ] (238.89,156.71) -- (217.86,141.42) -- (243.86,141.43) -- (264.89,156.72) -- cycle ;
\draw  [color=black ,draw opacity=1 ][fill=blue  ,fill opacity=0.35 ] (272.95,101.43) -- (264.93,76.7) -- (285.96,92) -- (293.98,116.73) -- cycle ;
\draw  [color=black ,draw opacity=1 ][fill=blue  ,fill opacity=0.35 ] (319.95,166.2) -- (311.93,141.46) -- (332.96,156.76) -- (340.98,181.49) -- cycle ;
\draw  [color=black ,draw opacity=1 ][fill=blue  ,fill opacity=0.35 ] (314.98,181.48) -- (293.95,166.18) -- (319.95,166.2) -- (340.98,181.49) -- cycle ;
\draw  [color=black ,draw opacity=1 ][fill=blue  ,fill opacity=0.35 ] (272.95,101.43) -- (264.91,126.16) -- (256.88,101.43) -- (264.93,76.7) -- cycle ;
\draw  [color=black ,draw opacity=1 ][fill=blue  ,fill opacity=0.35 ] (319.98,116.74) -- (293.98,116.73) -- (315.02,101.46) -- (341.02,101.47) -- cycle ;
\draw  [color=black ,draw opacity=1 ][fill=blue  ,fill opacity=0.35 ] (332.97,126.19) -- (311.93,141.46) -- (319.98,116.74) -- (341.02,101.47) -- cycle ;
\draw  [color=black ,draw opacity=1 ][fill=blue  ,fill opacity=0.35 ] (243.86,141.43) -- (217.86,141.42) -- (238.91,126.14) -- (264.91,126.16) -- cycle ;
\draw  [color=black ,draw opacity=1 ][fill=blue  ,fill opacity=0.35 ] (285.91,190.91) -- (264.86,206.18) -- (272.91,181.45) -- (293.95,166.18) -- cycle ;

\end{tikzpicture}
}
  \caption{Plain patch with 18 tiles}
  \label{fig:ring_plain}
\end{subfigure}%
\begin{subfigure}{.5\textwidth}
  \centering
  \scalebox{1.3}{
\begin{tikzpicture}[x=0.75pt,y=0.75pt,yscale=-1,xscale=1]

\draw  [color=black  ,draw opacity=1 ][fill=red  ,fill opacity=0.15 ] (388.05,116.78) -- (380,141.5) -- (354,141.49) -- (362.05,116.76) -- cycle ;
\draw  [color=black  ,draw opacity=1 ][fill=red  ,fill opacity=0.15 ] (354,141.49) -- (332.97,126.19) -- (341.02,101.47) -- (362.05,116.76) -- cycle ;
\draw  [color=black  ,draw opacity=1 ][fill=red  ,fill opacity=0.15 ] (340.98,181.49) -- (332.96,156.76) -- (354,141.49) -- (362.02,166.22) -- cycle ;
\draw  [color=black  ,draw opacity=1 ][fill=red  ,fill opacity=0.15 ] (388.02,166.23) -- (362.02,166.22) -- (354,141.49) -- (380,141.5) -- cycle ;
\draw  [color=black  ,draw opacity=1 ][fill=red  ,fill opacity=0.15 ] (293.98,116.73) -- (285.96,92) -- (307,76.73) -- (315.02,101.46) -- cycle ;
\draw  [color=black  ,draw opacity=1 ][fill=red  ,fill opacity=0.15 ] (264.89,156.72) -- (256.84,181.45) -- (230.84,181.43) -- (238.89,156.71) -- cycle ;
\draw  [color=black  ,draw opacity=1 ][fill=red  ,fill opacity=0.15 ] (264.91,126.16) -- (238.91,126.14) -- (230.88,101.41) -- (256.88,101.43) -- cycle ;
\draw  [color=black  ,draw opacity=1 ][fill=red  ,fill opacity=0.15 ] (306.93,206.2) -- (285.91,190.91) -- (293.95,166.18) -- (314.98,181.48) -- cycle ;

\draw  [color=black  ,draw opacity=1 ][fill=blue  ,fill opacity=0.15 ] (272.91,181.45) -- (264.86,206.18) -- (256.84,181.45) -- (264.89,156.72) -- cycle ;
\draw  [color=black  ,draw opacity=1 ][fill=blue  ,fill opacity=0.15 ] (238.89,156.71) -- (217.86,141.42) -- (243.86,141.43) -- (264.89,156.72) -- cycle ;
\draw  [color=black  ,draw opacity=1 ][fill=blue  ,fill opacity=0.15 ] (272.95,101.43) -- (264.93,76.7) -- (285.96,92) -- (293.98,116.73) -- cycle ;
\draw  [color=black  ,draw opacity=1 ][fill=blue  ,fill opacity=0.15 ] (319.95,166.2) -- (311.93,141.46) -- (332.96,156.76) -- (340.98,181.49) -- cycle ;
\draw  [color=black  ,draw opacity=1 ][fill=blue  ,fill opacity=0.15 ] (314.98,181.48) -- (293.95,166.18) -- (319.95,166.2) -- (340.98,181.49) -- cycle ;
\draw  [color=black  ,draw opacity=1 ][fill=blue  ,fill opacity=0.15 ] (272.95,101.43) -- (264.91,126.16) -- (256.88,101.43) -- (264.93,76.7) -- cycle ;
\draw  [color=black  ,draw opacity=1 ][fill=blue  ,fill opacity=0.15 ] (319.98,116.74) -- (293.98,116.73) -- (315.02,101.46) -- (341.02,101.47) -- cycle ;
\draw  [color=black  ,draw opacity=1 ][fill=blue  ,fill opacity=0.15 ] (332.97,126.19) -- (311.93,141.46) -- (319.98,116.74) -- (341.02,101.47) -- cycle ;
\draw  [color=black  ,draw opacity=1 ][fill=blue  ,fill opacity=0.15 ] (243.86,141.43) -- (217.86,141.42) -- (238.91,126.14) -- (264.91,126.16) -- cycle ;
\draw  [color=black  ,draw opacity=1 ][fill=blue  ,fill opacity=0.15 ] (285.91,190.91) -- (264.86,206.18) -- (272.91,181.45) -- (293.95,166.18) -- cycle ;

\node[circle,draw=red, fill=red, inner sep=1pt, minimum size= 4pt] (0) at (354,141.49) {};
\node[align=center] (p1) at (344,141.49) {\footnotesize $A$ };
\node[align=center] (p1) at (350,121) {\tiny $RI_5$ };
\node[align=center] (p1) at (350,161.49) {\tiny $RI_2$ };
\node[align=center] (p1) at (373,128.49) {\tiny $RI_4$ };
\node[align=center] (p1) at (373,154.49) {\tiny $RI_3$ };

\node[circle,draw=black, fill=black, inner sep=1pt, minimum size= 4pt] (0) at (319.98,116.74) {};
\node[align=center] (p1) at (310,125) {\footnotesize $B$ };
\node[align=center] (p1) at (326.98,123.74) {\tiny $S_5$ };
\node[align=center] (p1) at (317.98,108.74) {\tiny $S_2$ };

\node[circle,draw=black, fill=black, inner sep=1pt, minimum size= 4pt] (0) at  (272.95,101.43) {};
\node[align=center] (p1) at (278,114) {\footnotesize $C$ };
\node[align=center] (p1) at (282.95,96.43) {\tiny $S_1$ };
\node[align=center] (p1) at (265.95,101.43) {\tiny $S_3$ };

\node[circle,draw=black, fill=black, inner sep=1pt, minimum size= 4pt] (0) at (243.86,141.43) {};
\node[align=center] (p1) at (258.86,141.43) {\footnotesize $D$ };
\node[align=center] (p1) at (241.86,133.43) {\tiny $S_2$ };
\node[align=center] (p1) at (241.86,149.43) {\tiny $S_4$ };

\node[circle,draw=black, fill=black, inner sep=1pt, minimum size= 4pt] (0) at (272.91,181) {};
\node[align=center] (p1) at (278,170) {\footnotesize $E$ };
\node[align=center] (p1) at (279.91,188) {\tiny $S_5$ };
\node[align=center] (p1) at (265.95,181) {\tiny $S_3$ };

\node[circle,draw=black, fill=black, inner sep=1pt, minimum size= 4pt] (0) at  (319.95,166.2) {};
\node[align=center] (p1) at (311,158) {\footnotesize $F$ };
\node[align=center] (p1) at (317.95,174.2) {\tiny $S_4$ };
\node[align=center] (p1) at (329.95,161.2) {\tiny $S_1$ };

\node[circle,draw=black, fill=black, inner sep=1pt, minimum size= 4pt] (0) at  (307,76.73) {};
\node[align=center] (p1) at (320,76.73) {\footnotesize $G$ };
\node[align=center] (p1) at (303,96.73) {\tiny $RI_2$ };

\node[circle,draw=black, fill=black, inner sep=1pt, minimum size= 4pt] (0) at  (230.88,101.41) {};
\node[align=center] (p1) at (222,101.41) {\footnotesize $H$ };
\node[align=center] (p1) at (249.88,114.41) {\tiny $RI_3$ };

\node[circle,draw=black, fill=black, inner sep=1pt, minimum size= 4pt] (0) at  (230.84,181.43) {};
\node[align=center] (p1) at (222,181.43) {\footnotesize $I$ };
\node[align=center] (p1) at (249.84,168.43) {\tiny $RI_4$ };

\node[circle,draw=black, fill=black, inner sep=1pt, minimum size= 4pt] (0) at  (307,206.2) {};
\node[align=center] (p1) at (320,206.2) {\footnotesize $J$ };
\node[align=center] (p1) at (303,186.2){\tiny $RI_5$ };
\end{tikzpicture}
}
  \caption{Patch after the identification}
  \label{fig:ring_description}
\end{subfigure}
\caption{The original \textit{diamond ring} patch on the left and the same patch with indicated shallow holes (dots) and with a chosen ``origin" (red dot). The tiles are labelled with respect to the shallow hole they contain. The picture on the right shows the situation after step 3 of the algorithm. }
\label{fig:ring_labeling}
\end{figure}
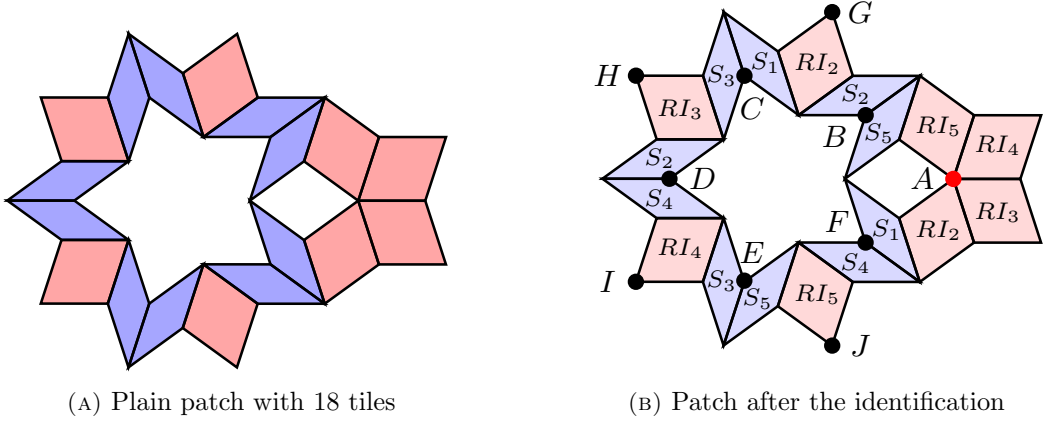

\begin{figure}[h!]
    \centering
    \includegraphics[width=0.4\textwidth]{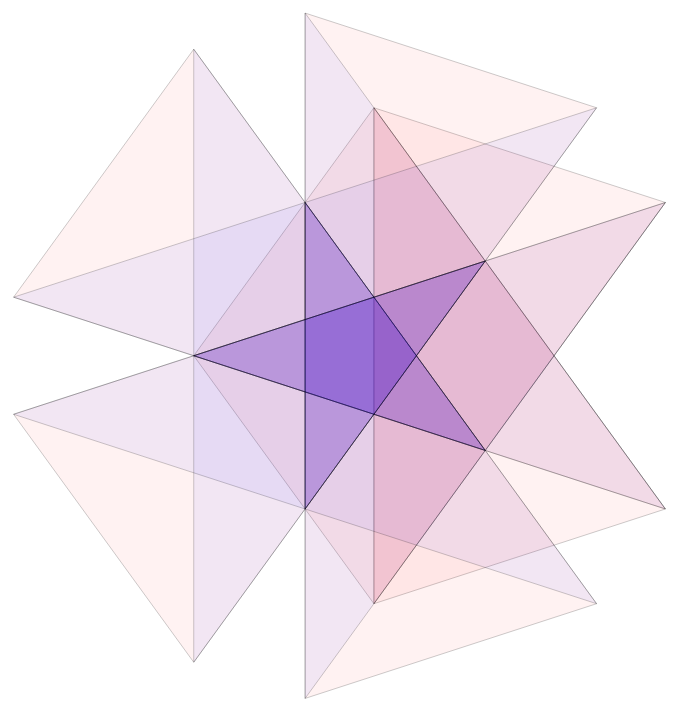}
    \caption{Intersection of dual tiles of the \textit{diamond ring} patch. They possess a~common intersection, the small violet triangle. }
    \label{fig:diam_intersection}
\end{figure}

Let us demonstrate the procedure on the following patch (this patch, called \emph{diamond ring}, supports an eigenfunction of a discrete Laplacian on the Penrose tiling, see \cite{DEFM22} for further details), see Figure \ref{fig:ring_plain}. This figure shows the initial data of the algorithm. Figure \ref{fig:ring_description} shows the result of step 2 (determining the shallow holes) and of step 3 (labelling the tiles). Table \ref{tab:rel_pos} summarises the paths from the `origin' (red point $A$) to (black) shallow holes (labelled with letters $B,\dots,J$), i.e., the relative translation vectors, i.e., the result of step 4. 
\begin{table}[h!]
\centering
\bgroup
\def\arraystretch{1.4}
\begin{tabular}{c|c}
     Shallow hole & Translation vector \\
     \hline \hline
     B & $\a_2 + \a_3 + \a_4$ \\
     C & $-\a_1 + \a_2 + 2\a_3+\a_4$\\
     D & $-\a_1 + 2\a_3 + 2\a_4$ \\
     E & $-\a_1 + \a_3 + 2\a_4 +\a_5$\\
     F & $\a_3 + \a_4 + \a_5$ \\
     G & $-\a_1 + \a_2 + \a_3 - \a_5 $ \\
     H & $ -2\a_1 + 2\a_3 + \a_4 - \a_5$ \\
     I & $-2\a_1 -\a_2 + \a_3 +2\a_4 $ \\
     J & $-\a_1 -\a_2 + \a_4 + \a_5 $ \\
     \hline\hline
\end{tabular}
\vspace{0.5cm}
\egroup
\caption{The positions of shallow holes of the \textit{diamond ring} patch relative to the ``origin" $A$. In particular, this is the result of step 4. For better readability, we abbreviate $\pi_{{}_\perp}(\bs{a}_i)$ to $\a_i$.  }
\label{tab:rel_pos}
\end{table}
Finally, Figure \ref{fig:diam_intersection} shows the result of the correspondence described in step 5, i.e., it depicts the corresponding dual triangles in the $\mathcal{S}^{\parallel}$ and their intersection (step 6), which is, in this particular case, a~triangle. Its area (step 7) is $\tfrac{1}{\sqrt{2(1165+521\sqrt{5})}}$. Thus, the frequency of the \textit{diamond ring} patch reads $\nu_{\mathrm{diam}} = \tfrac{1}{10}(34-21\tau) = \tfrac{1}{10}\tau^{-8} \in \mathcal{M}_{\mathscr{T}_{\mathrm{Pen}}}$. The total frequency of this patch (i.e., of all its possible rotates and space inversions) is $\nu^{\mathrm{tot}}_{\mathrm{diam}} = 34-21\tau$. 
We include other patches mentioned in \cite{DEFM22} in the Appendix.

The algorithm for obtaining patch frequencies can also be used for an entire class of tilings, namely, for those tilings obtained via the dualisation method. Usually, there is no need for distinguishing between deep and shallow holes, which makes the procedure slightly easier. On the other hand, another restriction may occur, but the idea and the basic scheme remain the same. By interchanging the roles of triangles and rhombuses, one can obtain the Tübingen Triangle Tiling (TTT) \cite{BKSZ90}. Using a different root lattice, one can also get patch frequencies for a plethora of quasiperiodic tilings with eight- and twelve-fold symmetry, including the Ammann--Beenker tiling \cite{BJ90,BJS91}. 

\section*{Acknowledgements}
I want to thank Michael Baake for introducing this problem to me, for valuable discussions and for all suggestions that helped to improve the manuscript. I would also like to thank Franz Gähler for explaining some properties of Ammann--Beenker tiling and to anonymous referees for several helpful comments. This work was supported by the German Research Foundation (DFG) within the CRC 1283/2 (2021 - 317210226) at Bielefeld University.

\newpage
\section*{Appendix 1 - Exact results for patches in Penrose tiling}
\begin{figure}[h!]
 \begin{subfigure}{.5\textwidth}
  \centering
  \includegraphics[width=0.4\textwidth]{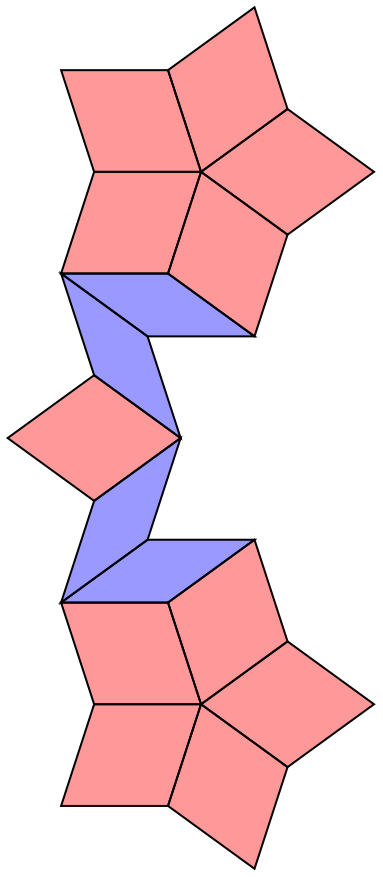}
\end{subfigure}%
\begin{subfigure}{.5\textwidth}
  \centering
  \includegraphics[width=0.8\textwidth]{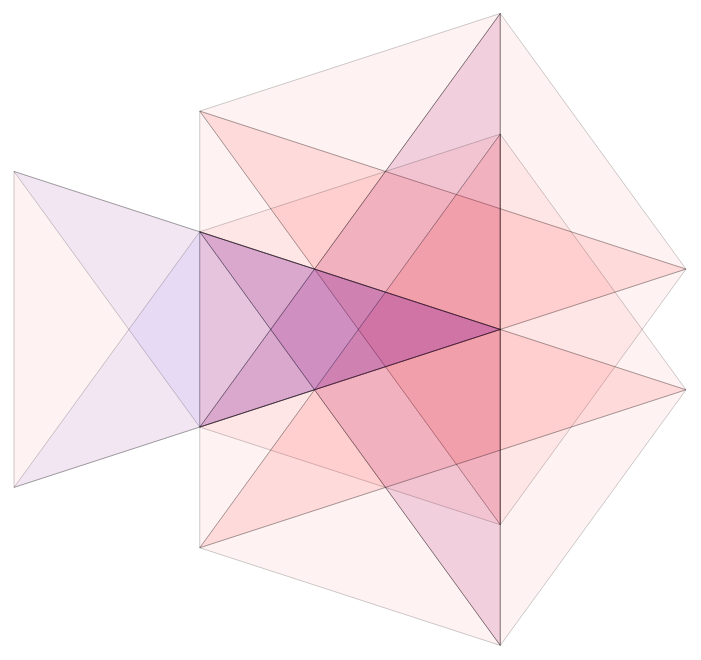}

\end{subfigure}
\caption{The \textit{two star} patch with 15 tiles. Its frequency is $\nu_{\mathrm{two}} = \tfrac{1}{5}(34-21\tau) = \tfrac{1}{5}\tau^{-8}$. }
\label{fig:two_star}
\end{figure}

\begin{figure}[h!]
 \begin{subfigure}{.5\textwidth}
  \centering
  \includegraphics[width=0.7\textwidth]{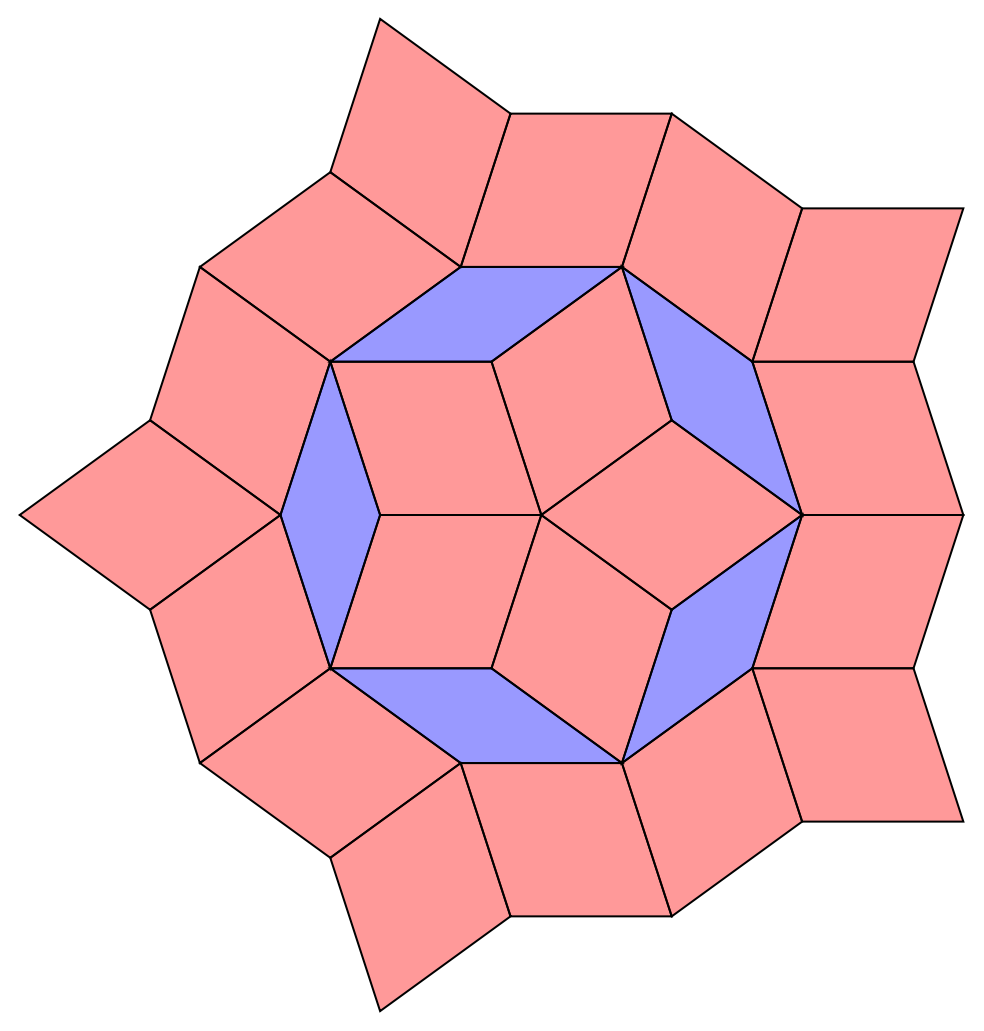}
\end{subfigure}%
\begin{subfigure}{.5\textwidth}
  \centering
  \includegraphics[width=0.8\textwidth]{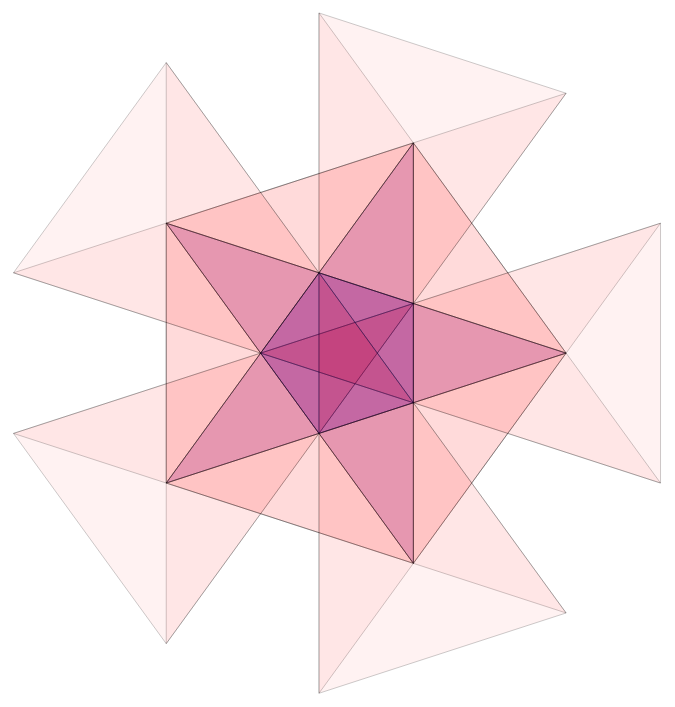}
\end{subfigure}
\caption{The \textit{filled circle} patch with 25 tiles. Its frequency is $\nu_{\mathrm{filled}} = \tfrac{1}{10}(123-76\tau) = \tfrac{2\tau-1}{10}\tau^{-9}$. }
\label{fig:filled circle}
\end{figure}

\begin{figure}[h!]
 \begin{subfigure}{.5\textwidth}
  \centering
  \includegraphics[width=0.7\textwidth]{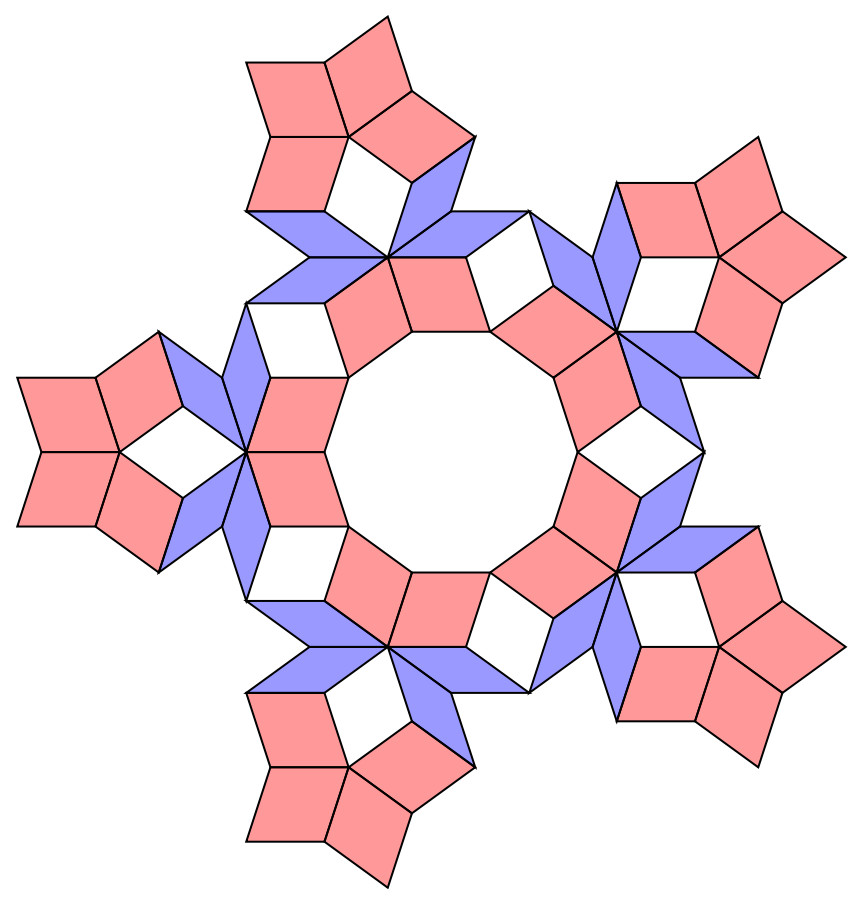}
\end{subfigure}%
\begin{subfigure}{.5\textwidth}
  \centering
  \includegraphics[width=0.77\textwidth]{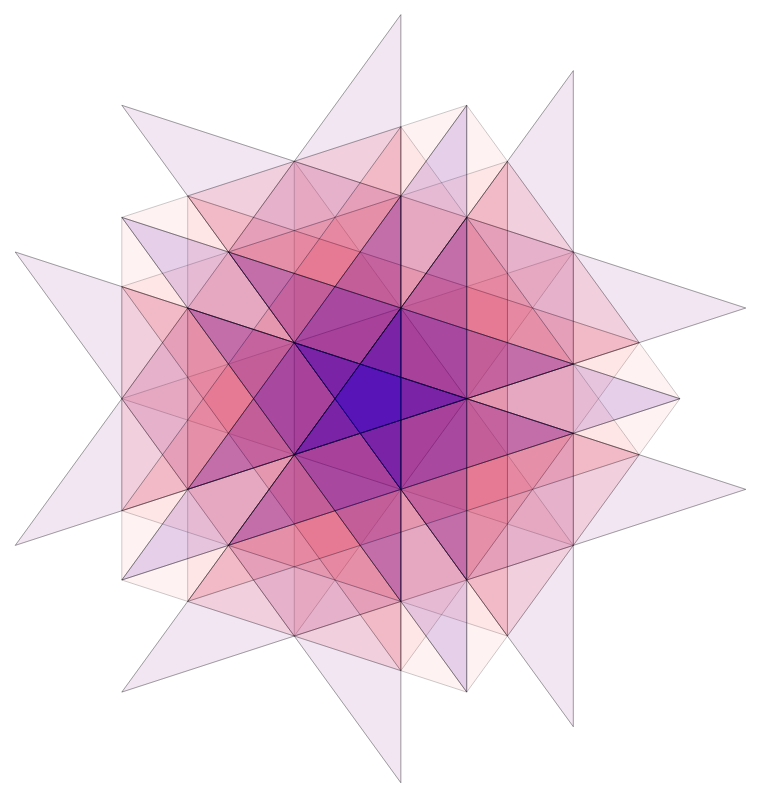}

\end{subfigure}
\caption{The \textit{big star} patch with 50 tiles. Its frequency is $\nu_{\mathrm{big}} = \tfrac{1}{10}(123-76\tau)=\tfrac{2\tau-1}{10}\tau^{-9}$. }
\label{fig:big_star}
\end{figure}

\begin{figure}[h!]
    \centering
    \includegraphics[width=0.72\textwidth]{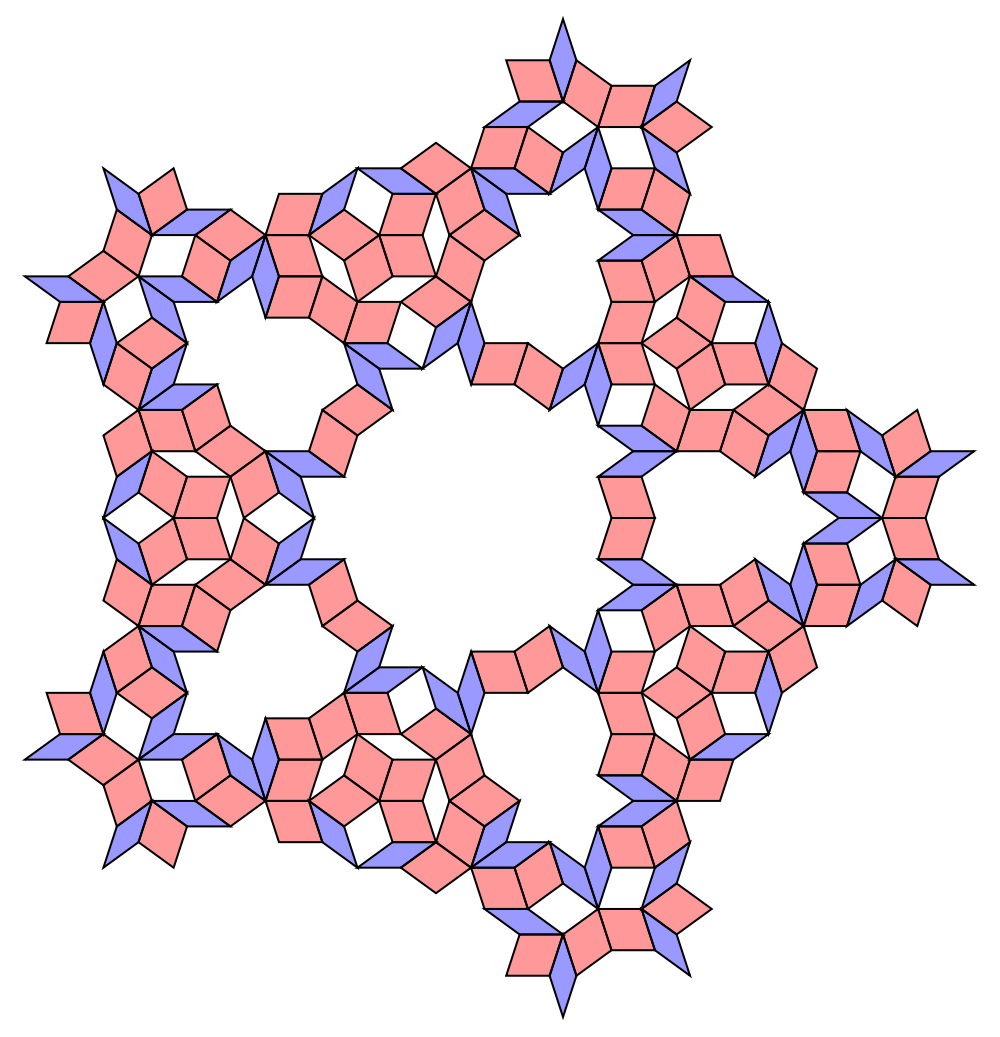}
    \caption{A 200-tiles patch}
    \label{fig:huge_1}
\end{figure}

\begin{figure}[h!]
    \centering
    \includegraphics[width=0.8\textwidth]{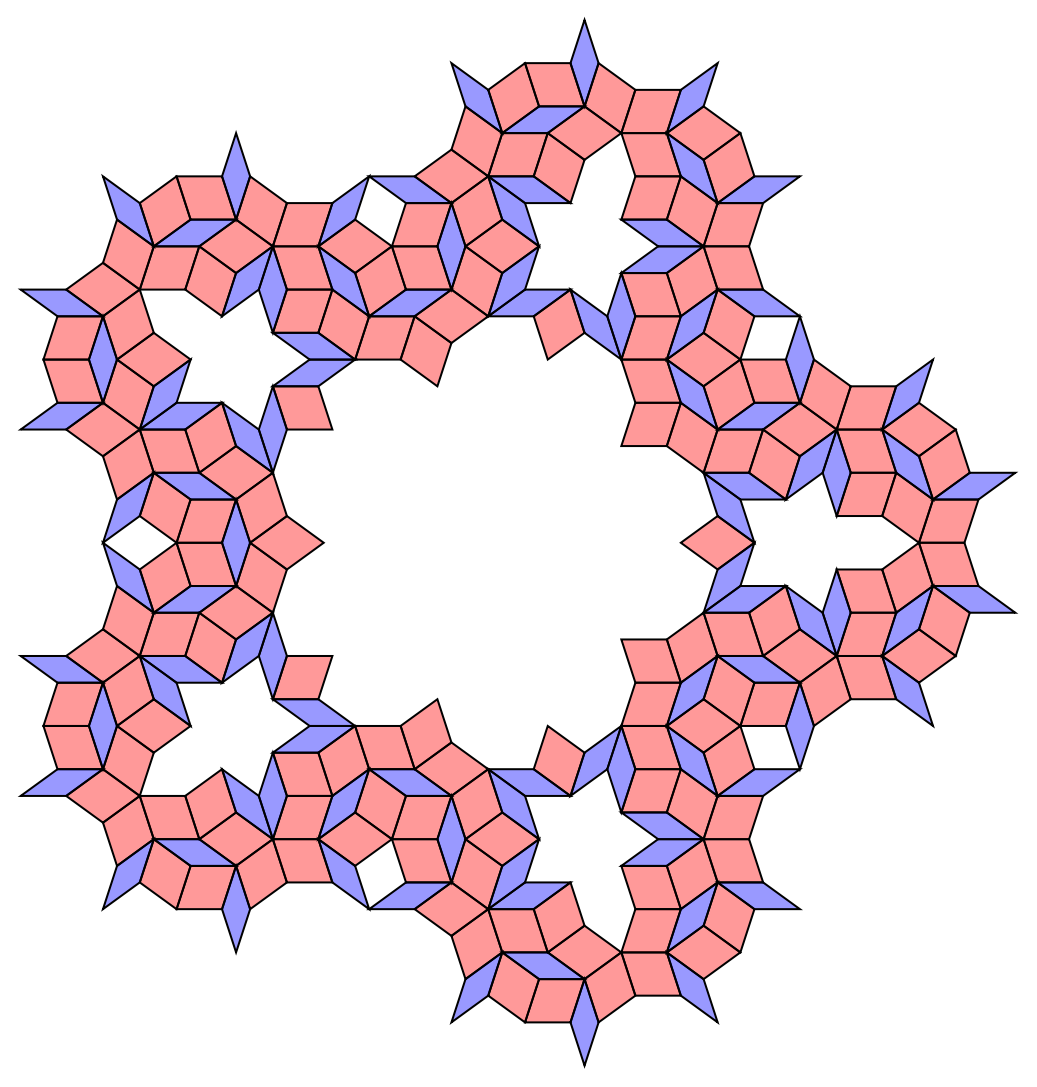}
    \caption{A 245-tiles patch}
    \label{fig:huge_2}
\end{figure}

In Figures \ref{fig:two_star}, \ref{fig:filled circle}, \ref{fig:big_star}, \ref{fig:huge_1}, and \ref{fig:huge_2}, we depict other patches that appear in \cite{DEFM22}, and the corresponding dual triangles in $\perp$-space. We also give the frequencies of these patches. 

\pagebreak

\begin{figure}[h!]
    \centering
    \includegraphics[width=0.95\textwidth]{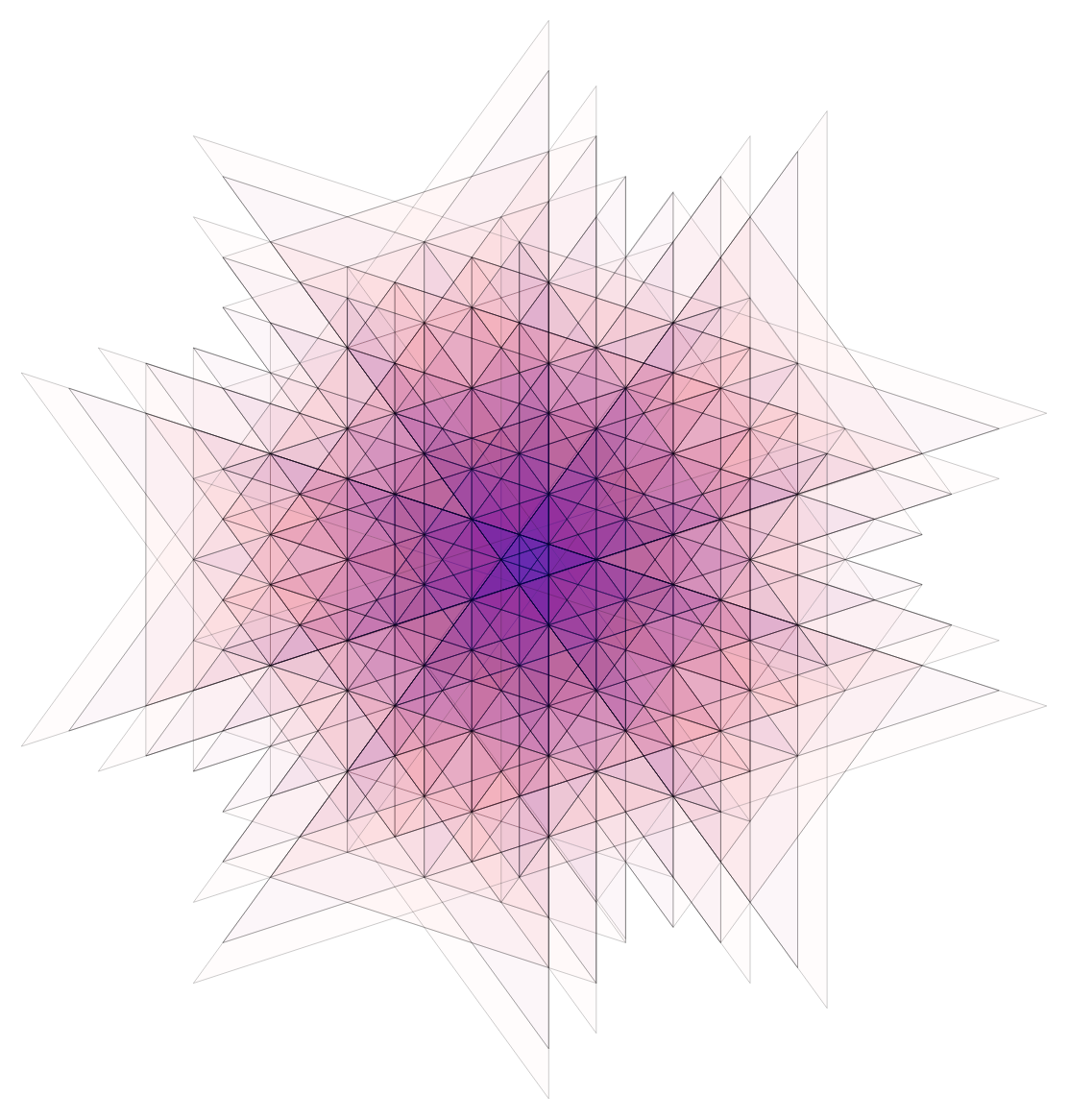}
    \caption{The dual image of the patch from Fig. \ref{fig:huge_1}. The frequency of this patch reads  $\nu_{\mathrm{huge1}} = \tfrac{1}{10}(2207-1364\tau) = \tfrac{2\tau-1}{10}\tau^{-15}$. }
    \label{fig:huge_1_window}
\end{figure}

\pagebreak

\pagebreak 

\begin{figure}[h!]
    \centering
    \includegraphics[width=1\textwidth]{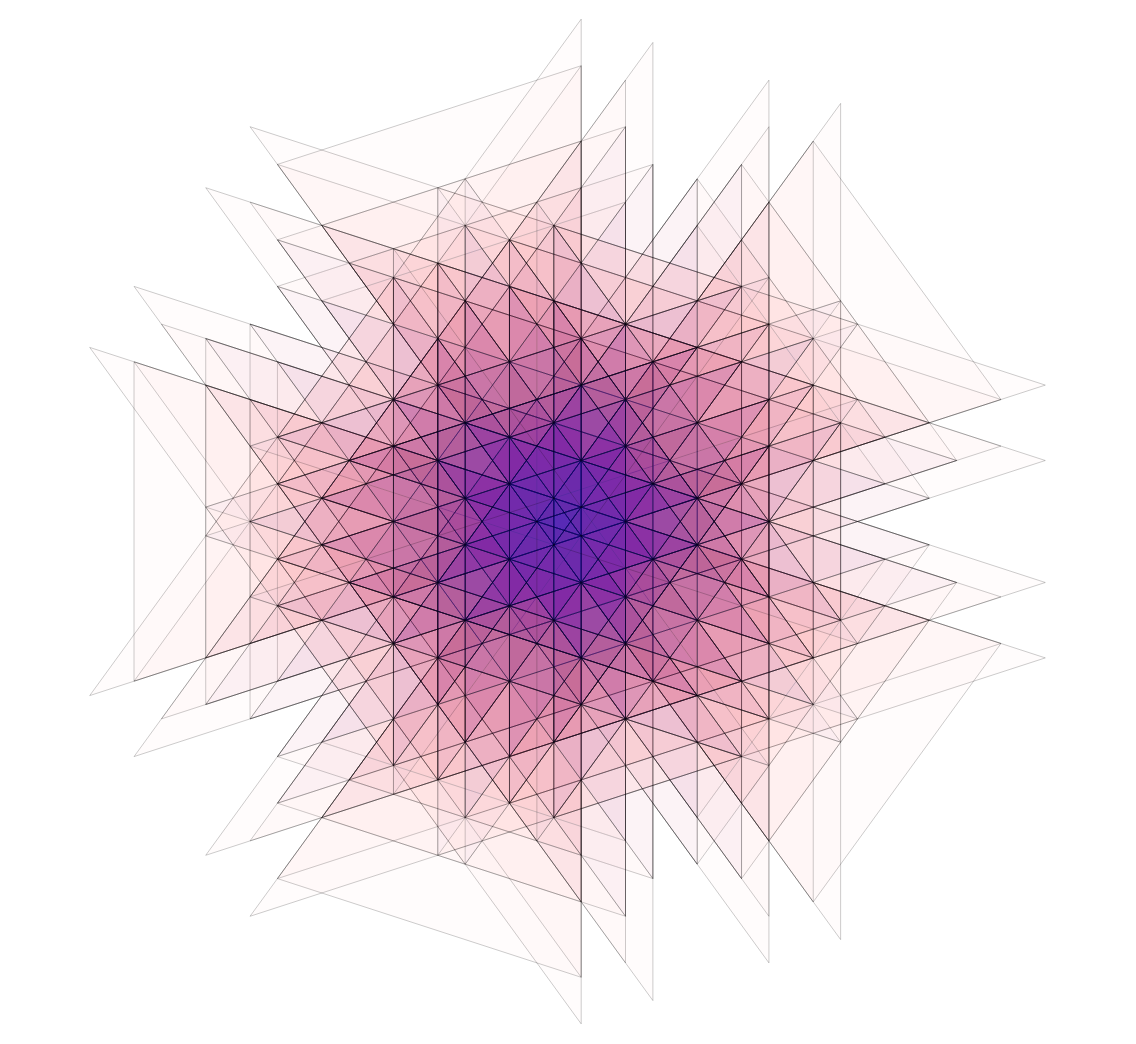}
    \caption{The dual image of the patch from Fig. \ref{fig:huge_2}. The frequency of this patch reads $\nu_{\mathrm{huge2}} = \tfrac{1}{10}(2207-1364\tau)= \tfrac{2\tau-1}{10}\tau^{-15}$. }
    \label{fig:huge_2_window}
\end{figure}

\pagebreak

\section*{Appendix 2 - Ammann--Beenker tiling}
Here, we briefly describe the setting for the Ammann--Bennker octagonal tiling. This tiling can be obtain via the dualisation of the 4-dimensional cubic lattice $\ZZ^4 = \langle \bs{e}_1,\bs{e}_2,\bs{e}_3,\bs{e}_4 \rangle_{\ZZ}$ which is self-dual. Recall that the Voronoi cell around the origin is the 4-cube given as 
\[ V_{\ZZ^4}(\bs{0}) \, = \, \bigg\{\bs{x}\in \RR^4 \ : \ | \bs{x}_i | \leqslant \myfrac{1}{2}    \mbox{ for all } \ 1\leqslant i \leqslant 4 \bigg \}. \]
The dual cells of the corresponding Voronoi complex  are of the form 
\[ V^{*}_{\ZZ^4}(\bs{0}^{*}) = \{ \bs{x}\in \RR^4 \ : \ 0\leqslant \bs{x}_i\leqslant 1   \mbox{ for all } \ 1\leqslant i \leqslant 4 \}  \, = \, V_{\ZZ^4}(\bs{0}) + \myfrac{1}{2}(1,1,1,1)^{T}.  \]

The symmetry group of the Voronoi cell is the \emph{hyperoctahedral group} $\Omega(4)$ \cite{BGO82}.  
All 2-boundaries of $V_{\ZZ^4}(\bs{0})$ are squares of the form
\[ Q(++\bigcirc\bigcirc) = \myfrac{1}{2}\{ \bs{e}_1+\bs{e}_2 +\lambda_1\bs{e}_3 + \lambda_2 \bs{e}_4 \ : \ -1\leqslant \lambda_1, \lambda_2 \leqslant 1  \}, \]
and all its possible images under the action of  $\Omega(4)$, which acts via permutations and sign flips. Together we obtain 24 congruent 2-boundaries. 
The dual boundaries are squares as well
\[ Q^{*}(++\bigcirc\bigcirc) = \{\mu_1 \bs{e}_1  + \mu_2\bs{e}_2 \ : \ 0\leqslant \mu_1, \mu_2 \leqslant 1 \}, \]
and the pairing of boundaries $Q$ and their dual boundaries $Q^*$ is one-to-one. 

\begin{figure}[h!]
    \centering
    \scalebox{0.9}{
\begin{tikzpicture}
			\filldraw[black] (0,0) circle (2pt) node[below left]{$\small\bs{0}$}; 
			\node[regular polygon, regular polygon sides=8, shape border rotate = -67.5,
			inner sep=1.4cm] (s) at (0,0) {} ;
			\draw[thick,  ->] (s.center) -- (s.corner 1)
			node [above]{$\small \pi_{{}_\parallel}(\bs{e}_1)$};
			\draw[thick, ->] (s.center) -- (s.corner 2)
			node [right]{ $\small \pi_{{}_\parallel}(\bs{e}_2)$};
			\draw[thick,  ->] (s.center) -- (s.corner 3)
			node [above right]{$\small \pi_{{}_\parallel}(\bs{e}_3)$};
			\draw[thick, ->] (s.center) -- (s.corner 4)
			node [above right]{$\small \pi_{{}_\parallel}(\bs{e}_4)$};
			\node (A) at (-1.5,-1.5) {$\parallel$-space};
			
			\filldraw[black] (8,0) circle (2pt) node[below left]{$\small\bs{0}$}; 
			\node[regular polygon, regular polygon sides=8, shape border rotate = -67.5,
			inner sep=1.4cm] (r) at (8,0) {} ;
			\draw[thick,  ->] (r.center) -- (r.corner 1)
			node [above]{$\small \pi_{{}_\perp}(\bs{e}_1)$};
			\draw[thick, ->] (r.center) -- (r.corner 4)
			node [right]{ $\small \pi_{{}_\perp}(\bs{e}_2)$};
			\draw[thick,  ->] (r.center) -- (r.corner 7)
			node [above right]{$\small \pi_{{}_\perp}(\bs{e}_3)$};
			\draw[thick, ->] (r.center) -- (r.corner 2)
			node [above right]{$\small \pi_{{}_\perp}(\bs{e}_4)$};
			\node (B) at (6.5,-1.5) {$\perp$-space};
\end{tikzpicture}
}
\caption{Projections of the standard basis $\bs{e}_1, \dots, \bs{e}_4$ into the two subspaces.  }
    \label{fig:basis_projections_AB}
\end{figure}
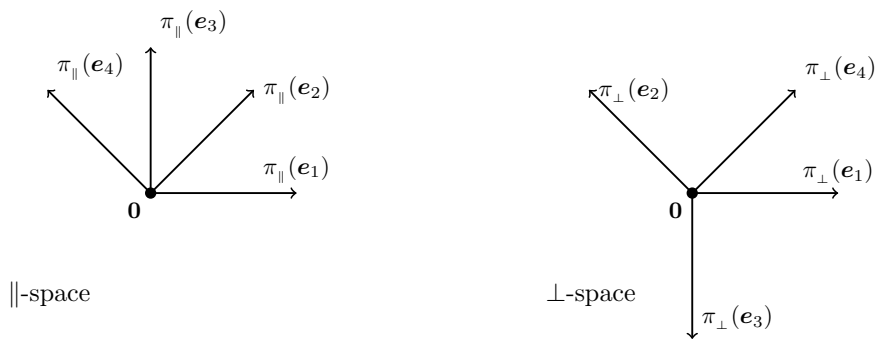

As in the case of the Penrose tiling, we need to find a suitable subgroup of the holohedry $\Omega(4)$ which possesses an (irreducible) representation in a plane. One can consider the dihedral group $D_8$ which is a proper subgroup of $\Omega(4)$.  This subgroup is generated by two elements $g_{_8}$, $s$ satisfying $g_{_8}^8 = s^2 = e$ and $(g_{_8}s)^2=e$. The generators act on the basis vectors $\bs{e}_i$ via the matrices
\[D(g_{_8}) = \begin{pmatrix} 0 & 0 & 0 & -1 \\ 1 & 0 & 0 & 0 \\ 0 & 1 & 0 & 0 \\ 0 & 0 & 1 & 0 \end{pmatrix}, \mbox{ and } D(s) = \begin{pmatrix} 1 & 0 & 0 & 0 \\ 0 & 0 & 0 & -1 \\ 0 & 0 & -1 & 0 \\ 0 & -1 & 0 & 0 \end{pmatrix}. \]
These matrices can be simultaneously brought to the real Jordan form, namely  
\[D(g_{_8}) \sim \begin{pmatrix} \cos \frac{\pi}{4} & -\sin \frac{\pi}{4} & 0 & 0 \\ \sin \frac{\pi}{4} & \cos \frac{\pi}{4} & 0 & 0 \\ 0 & 0 & \cos \frac{3\pi}{4} & -\sin \frac{3\pi}{4} \\ 0 & 0 & \sin \frac{3\pi}{4} & \cos \frac{3\pi}{4} \end{pmatrix}, \qquad D(s) \sim \begin{pmatrix} 1 & 0 & 0 & 0 \\ 0 & -1 & 0 & 0 \\ 0 & 0 & 1 & 0 \\ 0 & 0 & 0 & -1 \end{pmatrix}, \]
using the matrix 
\[ J \, = \, \sqrt{\myfrac{1}{2}} \begin{pmatrix} 
1 & \sqrt{\frac{1}{2}} & 0 & - \sqrt{\frac{1}{2}} \\
0 & \sqrt{\frac{1}{2}} & 1 &  \sqrt{\frac{1}{2}} \\
1 & -\sqrt{\frac{1}{2}} & 0 &  \sqrt{\frac{1}{2}} \\
0 & \sqrt{\frac{1}{2}} & -1 &  \sqrt{\frac{1}{2}} 
\end{pmatrix}. \] 
Taking the first two entries of each column of $J$, one gains the projections of the basis vectors into the $\parallel$-space, whereas taking the third and fourth one gives their $\perp$-projection. The projections are shown in Figure \ref{fig:basis_projections_AB}, and they already reveal the two shapes of tiles, namely a square, and a rhombus with the acute angle $\tfrac{\pi}{4}$.

\begin{figure}[h]
    \centering
    \scalebox{0.9}{
\tikzset{every picture/.style={line width=0.75pt}}        

\begin{tikzpicture}[x=0.75pt,y=0.75pt,yscale=-1,xscale=1]
%\draw  [color=red  ,draw opacity=0.5 ][line width=5.25]  (115.27,154.71) .. controls (122.17,161.61) and (125.66,170.63) .. (125.72,179.64) ;  
%\draw  [color=yellow  ,draw opacity=0.8 ][line width=5.25]  (320,144.74) .. controls (339.61,144.74) and (355.52,160.53) .. (355.52,180) ;  
\draw  [fill=yellow  ,fill opacity=0.35 ] (260.76,109.51) -- (190.36,179.91) -- (90.36,179.91) -- (161.1,109.15) -- (217.03,109.35) -- cycle ;

\draw  [draw opacity=0][line width=5.25]  (115.27,154.71) .. controls (122.17,161.61) and (125.66,170.63) .. (125.72,179.64) -- (90.33,179.64) -- cycle ;

\draw  [fill=red  ,fill opacity=0.35 ] (320,80) -- (420,80) -- (420,180) -- (320,180) -- cycle ;

\draw  [draw opacity=0][line width=5.25]  (320,144.74) .. controls (339.61,144.74) and (355.52,160.53) .. (355.52,180) -- (320,180) -- cycle ;

\draw    (380,180) -- (370,174) ;
\draw    (380,180) -- (370,186) ;
\draw    (380,80) -- (370,74) ;
\draw    (380,80) -- (370,86) ;
\draw    (320,120) -- (314,130) ;
\draw    (320,120) -- (326,130) ;
\draw    (420,120) -- (414,130) ;
\draw    (420,120) -- (426,130) ;
\draw    (130,180) -- (140,186) ;
\draw    (130,180) -- (140,174) ;
\draw    (220.45,109.64) -- (210.45,103.64) ;
\draw    (220.45,109.64) -- (210.45,115.64) ;
\draw    (232.54,137.61) -- (221.22,140.44) ;
\draw    (232.54,137.61) -- (229.71,148.92) ;
\draw    (118.65,151.62) -- (129.96,148.8) ;
\draw    (118.65,151.62) -- (121.48,140.31) ;

\end{tikzpicture}
}
    \caption{The decorated tiles for the Ammann--Beenker tiling. There are four different translation equivalent rhombus tiles and eight different square tiles. They differ by a rotation by an integer multiple of $\tfrac{\pi}{4}$. In the case of the rhombus tiles, one has to decide for suitable representatives since the decorated tile possesses a rotation symmetry by $\pi$. We decided to pick up as the representatives the rhombus on the picture, and its rotates by $\tfrac{\pi}{4}$, $\tfrac{\pi}{2}$, and $\tfrac{3\pi}{4}$.  }
    \label{fig:ab_tiles}
\end{figure}
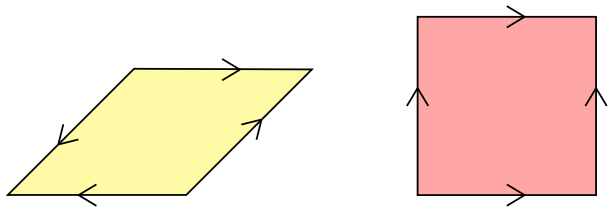

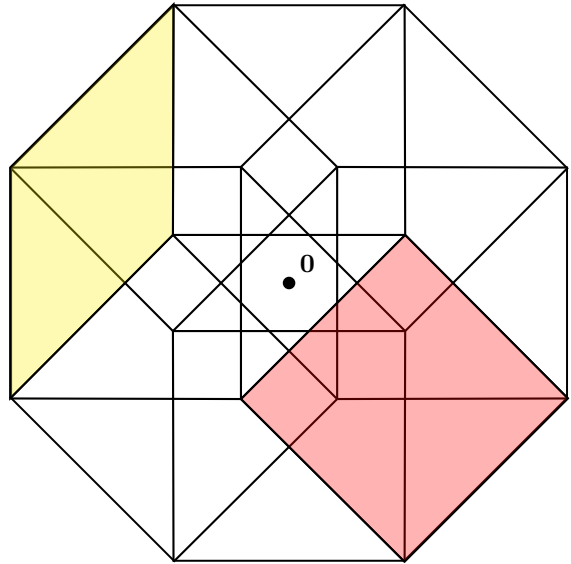
\begin{figure}[h!]
    \centering
    \scalebox{1}{
\begin{tikzpicture}
            \filldraw[black] (0,0) circle (2pt) node[above right]{$\small\bs{0}$};
			\node[regular polygon, regular polygon sides=8,
			inner sep=2.6131cm,draw] (S) at (0,0) {} ;
			\node[regular polygon, regular polygon sides=8,
			inner sep=1.0823cm] (s) at (0,0) {} ;
			\draw[thick] (S.corner 1) -- (s.corner 2) {};
			\draw[thick] (S.corner 2) -- (s.corner 3) {};
			\draw[thick] (S.corner 3) -- (s.corner 4) {};
			\draw[thick] (S.corner 4) -- (s.corner 5) {};
			\draw[thick] (S.corner 5) -- (s.corner 6) {};
			\draw[thick] (S.corner 6) -- (s.corner 7) {};
			\draw[thick] (S.corner 7) -- (s.corner 8) {};
			\draw[thick] (S.corner 8) -- (s.corner 1) {};
			\draw[thick] (S.corner 1) -- (s.corner 8) {};
			\draw[thick] (S.corner 2) -- (s.corner 1) {};
			\draw[thick] (S.corner 3) -- (s.corner 2) {};
			\draw[thick] (S.corner 4) -- (s.corner 3) {};
			\draw[thick] (S.corner 5) -- (s.corner 4) {};
			\draw[thick] (S.corner 6) -- (s.corner 5) {};
			\draw[thick] (S.corner 7) -- (s.corner 6) {};
			\draw[thick] (S.corner 8) -- (s.corner 7) {};
			\draw[thick] (s.corner 1) -- (s.corner 4) {};
			\draw[thick] (s.corner 2) -- (s.corner 5) {};
			\draw[thick] (s.corner 3) -- (s.corner 6) {};
			\draw[thick] (s.corner 4) -- (s.corner 7) {};
			\draw[thick] (s.corner 5) -- (s.corner 8) {};
			\draw[thick] (s.corner 6) -- (s.corner 1) {};
			\draw[thick] (s.corner 7) -- (s.corner 2) {};
			\draw[thick] (s.corner 8) -- (s.corner 3) {};
			\draw  [color=black  ,draw opacity=1 ][fill=yellow  ,fill opacity=0.3 ] (S.corner 2)  -- (S.corner 3)  -- (S.corner 4)  -- (s.corner 3)  -- cycle ;
			\draw  [color=black  ,draw opacity=1 ][fill=red  ,fill opacity=0.3 ] (S.corner 6)  -- (S.corner 7)  -- (s.corner 8)  -- (s.corner 5)  -- cycle ;
\end{tikzpicture}
}
\caption{Projection of the Voronoi cell $V_{\ZZ^4}(\bs{0})$ into the $\perp$-space with two 2-boundaries indicated. The yellow rhombus corresponds to $\pi_{{}_{\perp}}(Q(-+\bigcirc \bigcirc))$ and the red square is $\pi_{{}_{\perp}}(Q(+\bigcirc + \bigcirc))$. In contrast with the Penrose tiling, the centre of the window is placed in the origin. }
    \label{fig:AB_proj}
\end{figure}

The projections of the basis exhibit the desired octagonal symmetry. As above, we can project the 2-boundaries and get the Ammann--Beenker tiling as
\begin{equation}
\label{eq:tiling_descr_AB}
\mathscr{T}_{\mathrm{AB}}(\bs{c}_{{}_\perp}) \, = \, \left\{ \pi_{{}_\parallel}(Q^{*}) \ : \ Q^{*} \mbox{ is a 2-boundary},  \ \bs{c}_{{}_\perp} \in \pi_{{}_\perp} (Q) \right\}. 
\end{equation}
We choose the vector $\bs{c}_{{}_\perp}$ so that it does not belong to any 1-boundary of any Voronoi cell, similarly to the Penrose case. 
In contrast to the Penrose tiling, we project the dual boundaries into the $\parallel$-space, but this does not cause any difficulties. The $\perp$-projection of the Vorornoi cell with projections of two particular 2-boundaries is shown in Figure \ref{fig:AB_proj}. The area of the projection (which is an octagon) is $1+\sqrt{2}$. Up to a translation, we have twelve different tiles --- four rhombuses and eight squares (!). This, perhaps surprising, fact follows from the decorations of the Ammann--Beenker tiles, see \cite{TAO} for further details. Figure \ref{fig:ab_tiles} shows the tiles and their decorations. 

\begin{figure}
    \centering
    \input{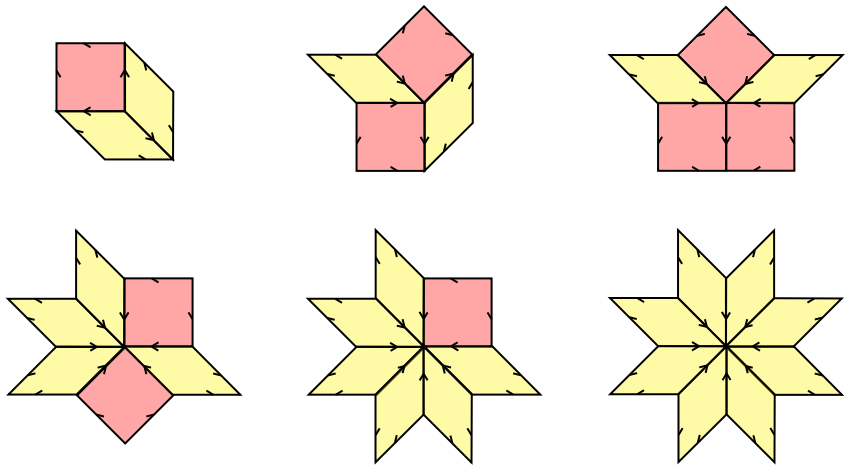}
    \caption{All allowed vertex configurations (up to rotations) within THE Ammann--Beenker tiling displayed with decorations. }
    \label{fig:ab_config}
\end{figure}

As in the case of the Penrose tiling, we can determine all elementary polygons and obtain all possible vertex configurations as shown in Figure \ref{fig:ab_config}.

Since there are no holes in this setting (as $\ZZ^4$ is self-dual as a lattice), the algorithm for determining the patch frequencies has to be modified as follows. One has to replace `the distinguishing between deep and shallow holes' in step 1 by `decorating the tiles', and in step 3, one has to replace `any shallow hole' with `any vertex point', since there is only a single translation class.
And, of course, in the last step, one has to divide by the accurate area of the window, in this case by $1+\sqrt{2}$.  
No other changes are needed. The patch frequencies are contained in the frequency module  $\mathcal{M}_{\mathscr{T}_{\mathrm{AB}}}$ which reads $\mathcal{M}_{\mathscr{T}_{\mathrm{AB}}} = \tfrac{1}{2} \ZZ[\lambda]$ with $\lambda = 1+\sqrt{2}$, the \emph{silver mean} \cite[Ex. 7.9]{TAO}. 

We enclose several patches of the Ammann--Beenker tiling which appear in \cite{DEFM22} with their frequencies. 

\begin{figure}[h!]
 \begin{subfigure}{.5\textwidth}
  \centering
  \includegraphics[width=0.7\textwidth]{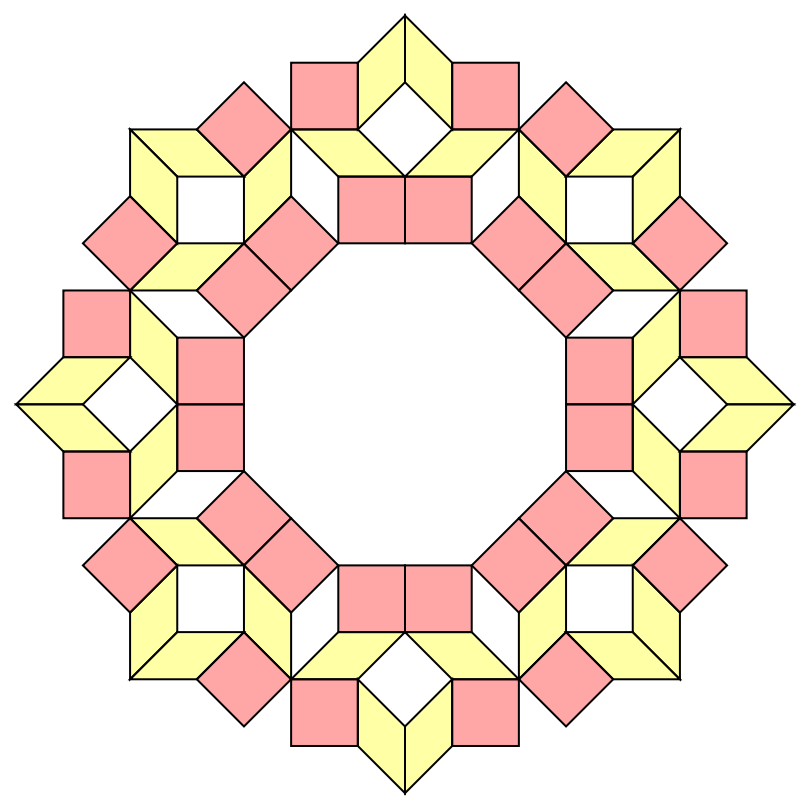}
\end{subfigure}%
\begin{subfigure}{.5\textwidth}
  \centering
  \includegraphics[width=0.7\textwidth]{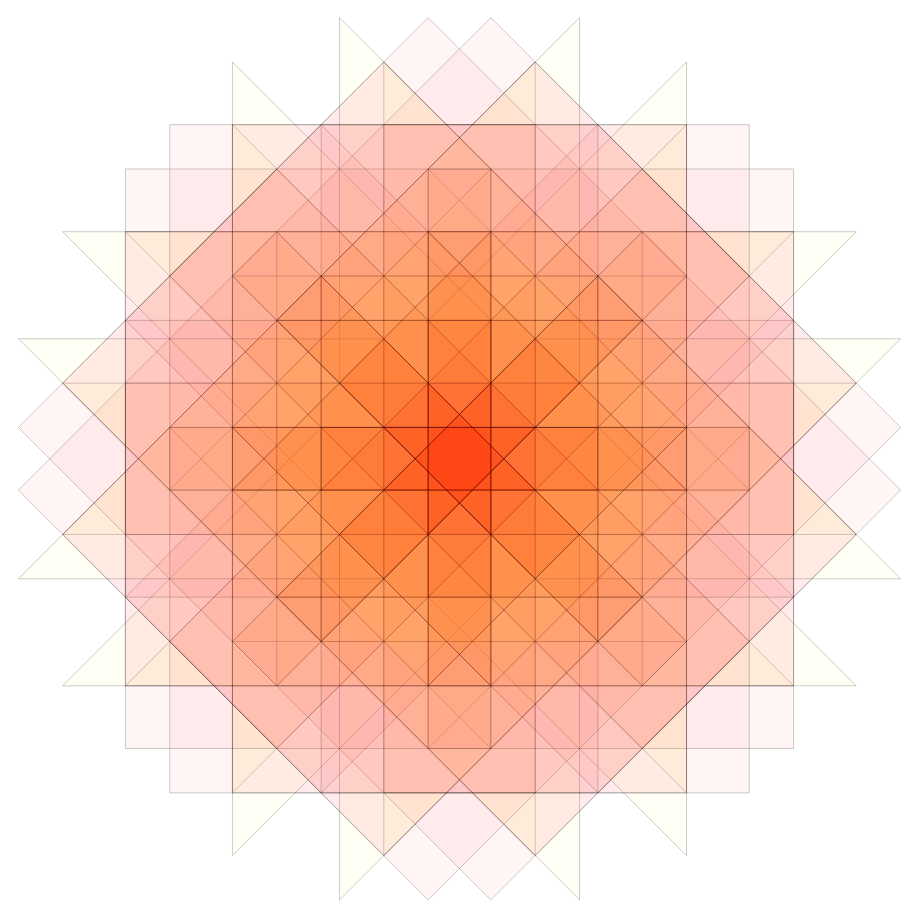}

\end{subfigure}
\caption{A 64-tiles patch. Its frequency is $\nu_{\mathrm{64}} = 29\lambda -70= \lambda^{-5}$. }
\label{fig:ab_64}
\end{figure}

\begin{figure}[h!]
    \centering
    \includegraphics[width=0.58\textwidth]{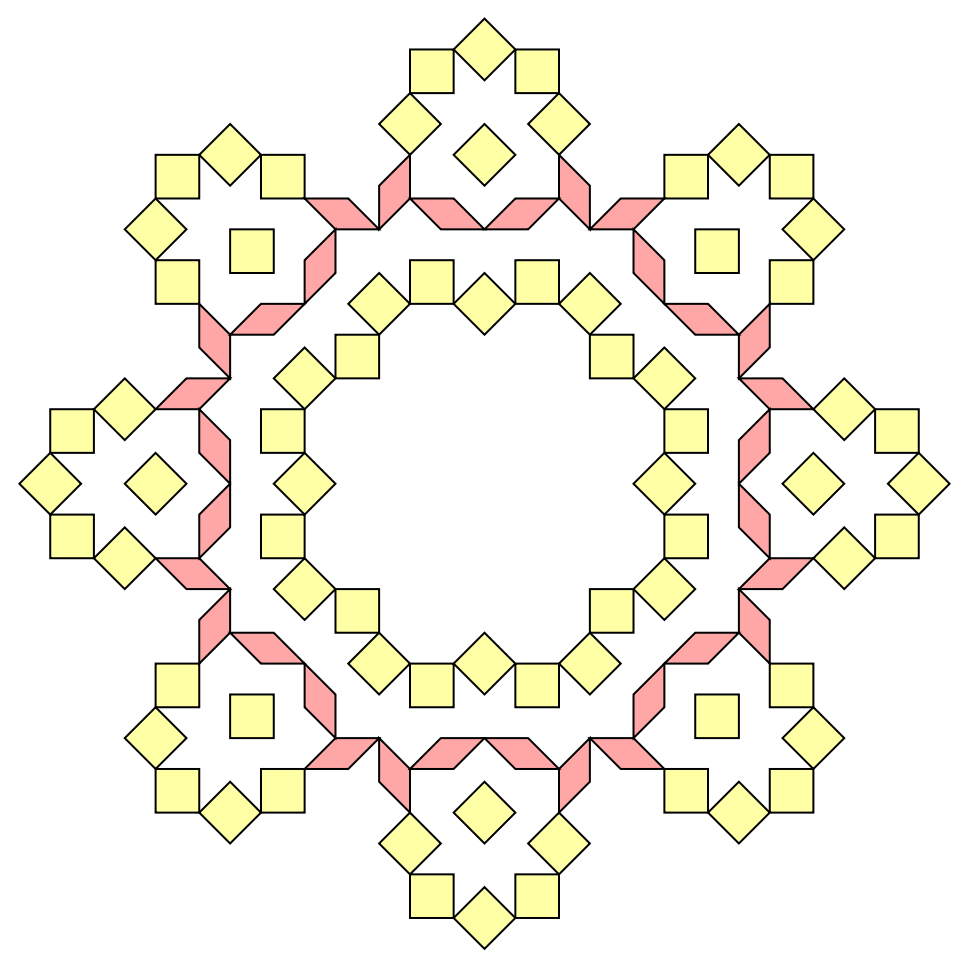}
    \caption{A 104-tiles patch. Its frequency is $\nu_{\mathrm{104}} = 985-408\lambda = \lambda^{-8}$. }
    \label{fig:ab_104}
\end{figure}

\begin{figure}[h!]
    \centering
    \includegraphics[width=1\textwidth]{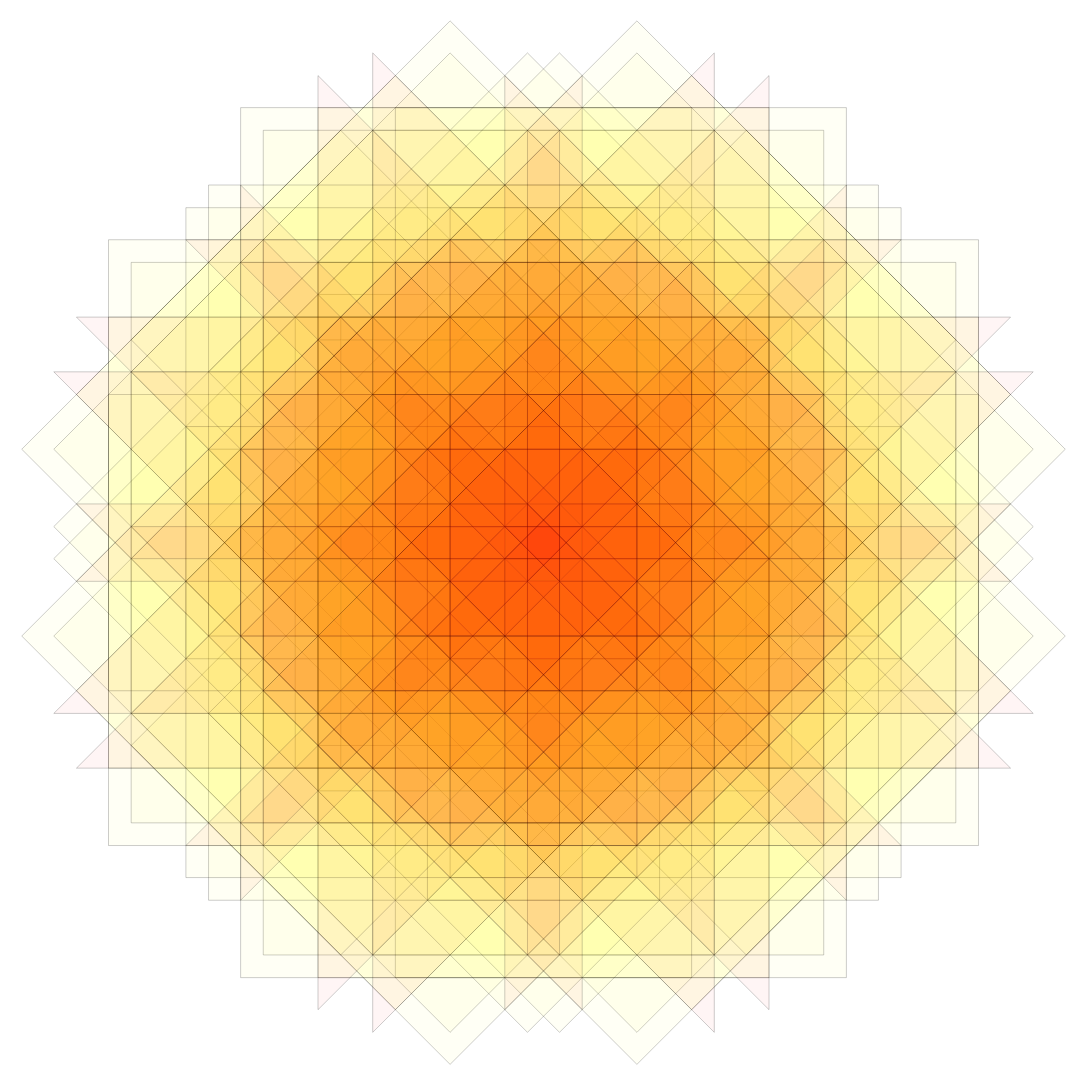}
    \caption{The intersection of dual tiles of the patch from Figure \ref{fig:ab_104}.}
    \label{fig:ab_104_win}
\end{figure}

\begin{figure}[h!]
    \centering
    \includegraphics[width=0.95\textwidth]{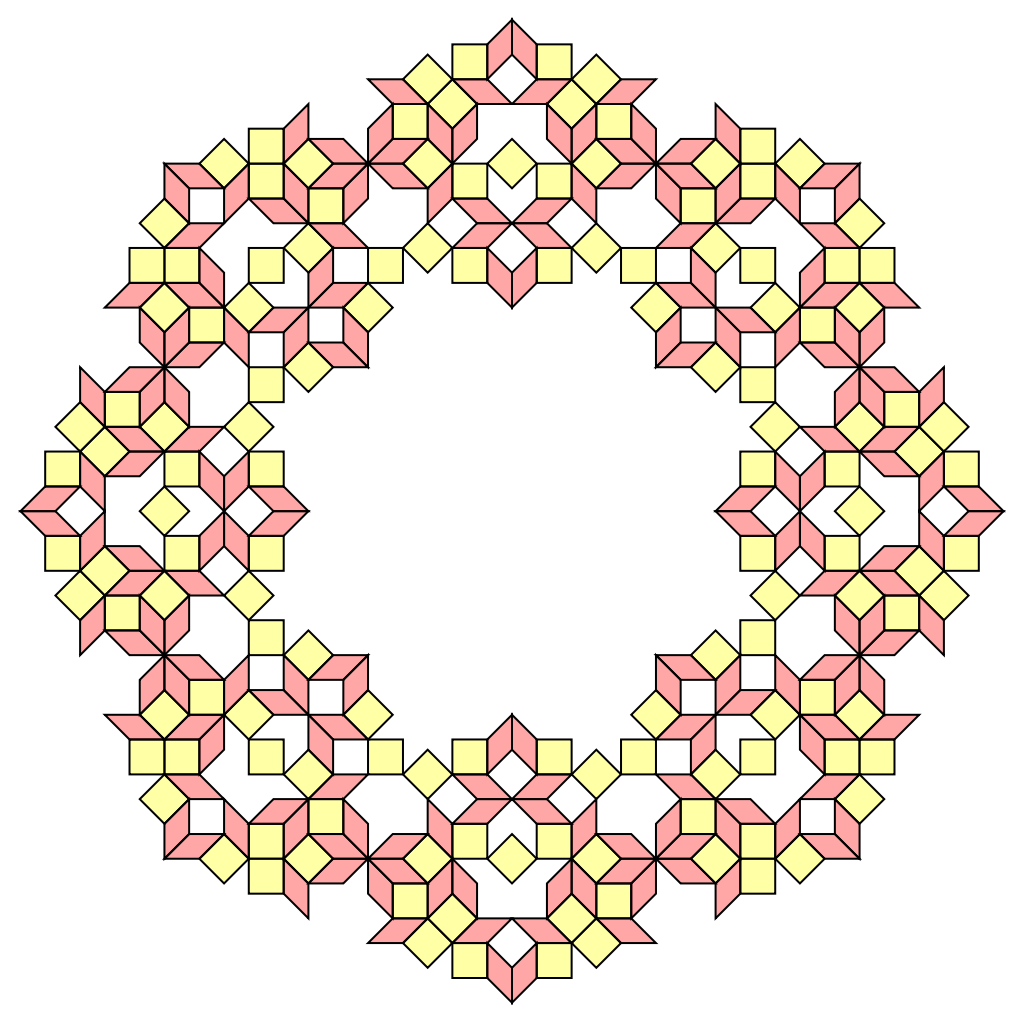}
    \caption{A 328-tiles patch. Its frequency is $\nu_{\mathrm{328}} = 985-408\lambda= \lambda^{-8}$. }
    \label{fig:ab_328}
\end{figure}

\begin{figure}[h!]
    \centering
    \includegraphics[width=0.95\textwidth]{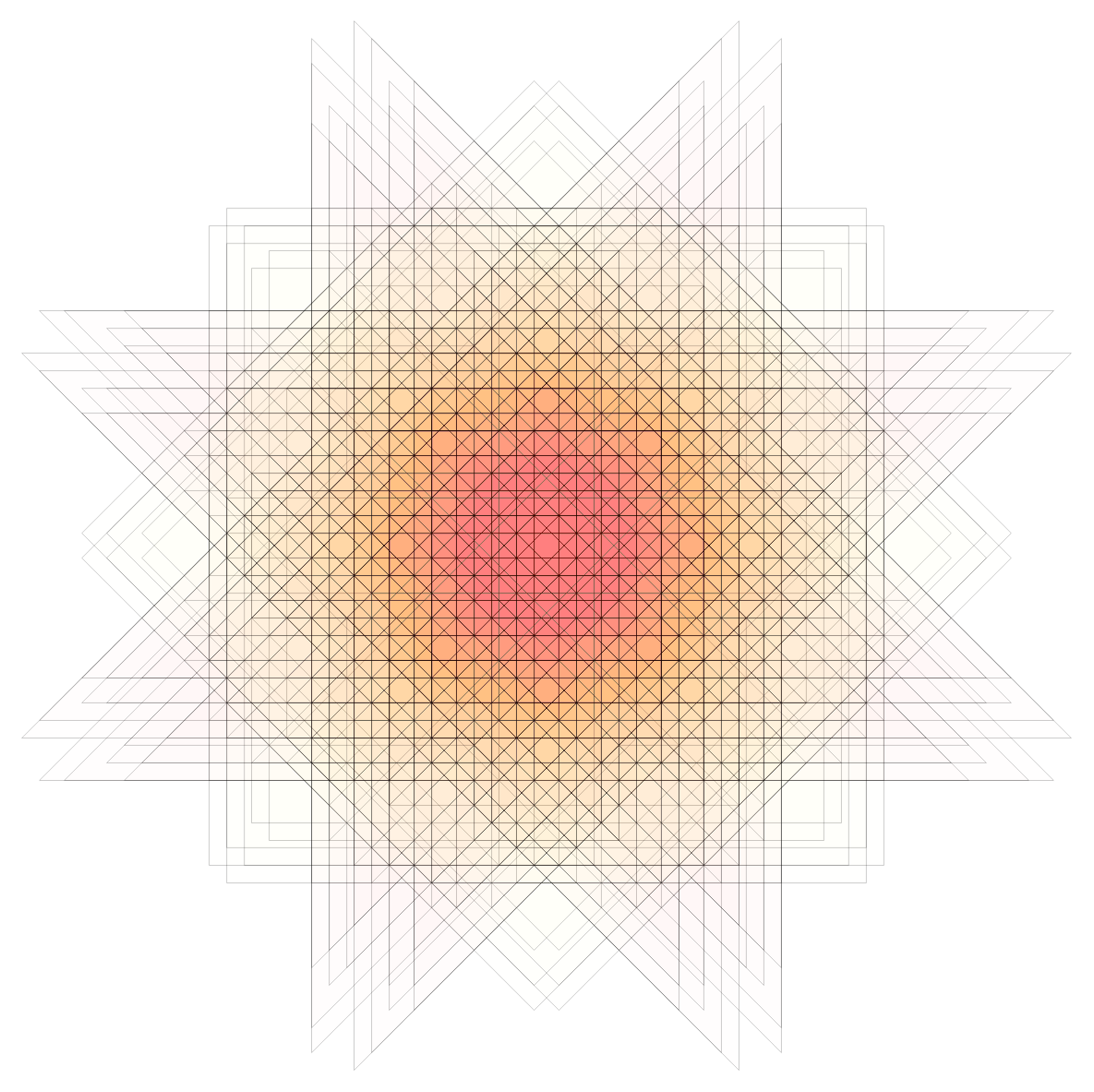}
    \caption{The intersection of dual tiles of the patch from Figure \ref{fig:ab_328}.}
    \label{fig:ab_328_win}
\end{figure}

\end{document}